\documentclass[preprint,12pt]{elsarticle}

\usepackage{amssymb}
\usepackage{amsmath}

\usepackage{amssymb,amsmath,amsfonts,latexsym,stackrel} 
\usepackage{hyperref}

 \newtheorem{thm}{Theorem}
 \newproof{pf}{Proof}

\journal{}

\begin{document}

\begin{frontmatter}



\title{Sibuya probability distributions
and numerical evaluation of fractional-order operators}


\author[1]{Nikolai Leonenko\corref{cor1}%
               \fnref{fn1}\fnref{fn2}} 
\ead{LeonenkoN@cardiff.ac.uk}

\affiliation{organization={Cardiff University, School of Mathematics},
            addressline={Abacws, Senghennydd Road}, 
            city={Cardiff},
            postcode={CF~24~4AG}, 
            country={United Kingdom}}

\author[2]{Igor Podlubny\fnref{fn1}\fnref{fn3}} 
\ead{igor.podlubny@tuke.sk} 

\affiliation{organization={Technical University of Kosice, BERG Faculty},
            addressline={Nemcovej 3}, 
            city={Kosice},
            postcode={04200}, 
            country={Slovakia}}

\cortext[cor1]{Corresponding author}
\fntext[fn1]{Nikolai Leonenko (NL) and Igor Podlubny (IP) would like to thank for support and hospitality during the programmes ``Fractional Differential Equations'',  ``Uncertainly Quantification and Modelling of Materials'', both supported by EPSRC grant EP/R014604/1, and the programme ``Stochastic systems for anomalous diffusion'', supported by EPSRC grant EP/Z000580/1, at Isaac Newton Institute for Mathematical Sciences, Cambridge. } 
\fntext[fn2]{Also NL was partially supported under the ARC Discovery Grant DP220101680 (Australia), Croatian Scientific Foundation (HRZZ) grant ``Scaling in Stochastic Models'' (IP-2022-10-8081), grant FAPESP 22/09201-8 (Brazil) and the Taith Research Mobility grant (Wales, Cardiff University). Also, NL would like to
thank University of Rome ``La Sapienza'' for hospitality as Visiting Professor (June
2024) where the paper was initiated.}
\fntext[fn3]{The work of IP is also partially supported by grants VEGA 1/0674/23, APVV-22-0508, and ARO W911NF-22-1-0264.}

\begin{abstract}
In this work we explore the Sibuya discrete probability distribution, which serves as the basis and the main instrument 
for numerical simulations of Gr\"unwald--Letnikov fractional derivatives by the Monte Carlo method.
We provide three methods for simulating the Sibuya distribution.
We also introduce the Sibuya-like sieved probability distributions, and apply them to numerical 
fractional-order differentiation. 
	Additionally, we use the Monte Carlo method for evaluating fractional-order integrals,
and suggest the notion of the continuous Sibuya probability distribution. 
	The developed methods and tools are illustrated by examples of computation.	
We provide the MATLAB toolboxes for simulation of the Sibuya probability distribution, 
and for the numerical examples. 
\end{abstract}

\begin{keyword}
fractional calculus (primary) \sep
fractional differentiation \sep numerical computations \sep Monte Carlo method \sep stochastic processes


\MSC 26A33 (primary) \sep 65C05 \sep 65D25

\end{keyword}

\end{frontmatter}



\section{Introduction} \label{sec:1}

A class of discrete time random walks has been recently introduced for numerical solution of fractional-order partial differential equations, including fractional-order subdiffusion equations with time fractional derivatives of order $\alpha \in (0, 1)$ (see \cite{Nichols-2018,DH-1,DH-2} and the references therein).  It turns out that the main role in Monte Carlo method for subdiffusion is played by the Sibuya distribution. 
	However, some problems require the order of the time fractional derivative is greater than one, and thus such problems require the Monte Carlo method for fractional differential operators of order $\alpha > 1$. For example, in fractional-order diffusion-wave equation the order is $ \alpha \in (1, 2]$. In such case, the signed probability distributions appear, as in \cite{LP-MCFD-2}. 

This paper is a continuation of \cite{LP-MCFD-1,LP-MCFD-2} 
and is devoted to using signed discrete distributions 
(or extended, or negative probabilities -- the terminology has not been unified yet) 
in the framework of the Monte Carlo method for numerical evaluation of fractional-order derivatives.  

Our framework is close to the work of Zhang, Li, and Kerns (see  Corollary 3.1 (3$^{o}$) in \cite{ZLK-2017}), and also to the one studied by Burgin \cite{Burgin-2012,Burgin-2013}. 

Burgin studies a static environment, where the state space $\Omega$
can be divided in two irreducible parts $\Omega_{+}$ and $\Omega_{-}$. 
The elements of $\Omega_{-}$ are called anti-events, and they are usually 
connected to negative objects: encounting a negative object is a negative event.
The author then states that an example of negative objects is given 
by antiparticles, the antimatter counterpart of quantum particles. 
Anti-events are assigned negative probabilities. 

Another interesting interpretation is given by Sz\'ekely~\cite{Szekely-2005},
who proves that we can encounter a negative probability if we work
with random variable having a signed distribution. 
In addition, if $X$ has a signed distribution, 
then there exist two random variables $Y$ and $Z$ 
having ordinary (not signed) such that $X+Y=Z$ in distribution.
Therefore, $X$ can be seen as a ``difference'' of two ordinary 
random variables $Z$ and $Y$. 

Signed distributions are related to the concept of subjective probability theory
given by De~Finetti \cite{DeFinetti-2017}. 

In quantum physics Wigner \cite{Wigner-1932}, Dirac \cite{Dirac-1942}, 
and Feynman \cite{Feynman-1987} allow probabilities (and also energies) to be negative,
in order to side-step the uncertainty principle of quantum mechanics.

Concepts of virtual events are known also in psychology \cite{Bollen-2002},
economics \cite{Hu-2017}, and medicine \cite{Rabe-Hesketh-Skrondal-2008}.
Some additional references can be found in \cite{Caprio-Mukherjee-arxiv-2023}.

The paper is organized as follows. 

We start with recalling how the Sibuya distribution is related to fractional-order differentiation, 
and focus on developing the methods for simulation of the Sibuya distribution for $0 < \alpha < 1$. 
We derive three different methods 
and provide their implementation in a freely downloadable MATLAB toolbox;  
all three methods are perfectly comparable, and any of them can be used for practical simulation. 

Then we focus on fractional-order differentiation of orders higher than one. 
In this case the Sibuya distribution becomes a signed probability distribution,
and we develop a general Monte Carlo method to numerical evaluation of 
fractional-order derivatives of order higher than one 
using separation of positive and negative probabilities. 

After that we turn our attention to fractional-order integration. 
It appears that fractional-order integrals can also be evaluated 
by the Monte Carlo method using a continuous probability distribution 
that can be called the continuous Sibuya distribution. 

The developed methods are illustrated by a set of examples,  
in which the results of computations using the Monte Carlo method are compared 
with the exact fractional-order derivatives and integrals of the corresponding functions. 
The examples are also provided as a MATLAB toolbox.

\section{Gr\"unwald--Letnikov fractional derivative} \label{sec:2}

Let us recall that the Gr\"unwald--Letnikov fractional derivative is a non-local operator   defined as \cite{Podlubny-FDE} 
\begin{equation} \label{GL-definition}
	_{0}D_{t}^\alpha f(t) = \lim_{h \rightarrow 0} A_{h}^{\alpha} f(t), 
	\quad 
	\alpha > 0, \quad h>0,
\end{equation}
whenever the limit exists, where
\begin{equation}
	A_h^\alpha f(t) = \frac{1}{h^\alpha} \left( \Delta_h^\alpha f \right) (t)
\end{equation}
and
\begin{equation} 
	\left( \Delta_h^\alpha f \right) (t) = 
	\sum_{k=0}^{\infty}
	(-1)^k  {\alpha \choose k}	
	f(t-kh),
\end{equation}
for sufficiently good functions $f$ \cite{Podlubny-FDE}.

\noindent
Here appear the binomial coefficients
\begin{equation}
	{\alpha \choose k} = 
	\frac{(\alpha)_k}{k!},
\end{equation}
where 
\begin{equation}
	(\alpha)_k  = 
	\frac{\Gamma (\alpha+1)}{\Gamma (\alpha+1-k)} = 
	\alpha (\alpha-1) \ldots (\alpha - k + 1)
\end{equation}
is the Pochhammer's symbol.

It is known that for every bounded function $f$ the semigroup property holds, i.e.
\begin{equation}
	\Delta_h^\alpha 
	\Delta_h^\beta  f
	= 
	\Delta_h^{\alpha+\beta } f,
	\quad
	(\alpha, \beta > 0),
\end{equation}
and the Gr\"unwald--Letnikov derivative in many cases coincides with 
the Riemann--Liouville, Caputo, Weyl or Marchaud fractional derivatives 
\cite{Podlubny-FDE,Ferrari-2018,LP-MCFD-1}.

\section{Sibuya distribution} \label{sec:3}

Consider the sequence of independent Bernoulli trials, in which the $k$-th trial 
has probability of success $\alpha/k$ ($0<\alpha<1$, \, $k = 1, 2, \ldots$). 
Let $Y \in \{1, 2, \ldots \}$ be the trial number in which the first success occurs.
Then 
$$
	\mathbb{P} (Y=1) = p_1 = p_1(\alpha) = \alpha,
$$
$$
	 \mathbb{P} (Y=k) = p_k = p_k(\alpha) =
	 (1-\alpha) 
	 (1-\frac{\alpha}{2})
	\ldots
	(1-\frac{\alpha}{k-1}) 
	\frac{\alpha}{k},
	\quad k \geq 2,
$$	
or, in other words, its probability mass function is
\begin{eqnarray}
	p_k = \mathbb{P} (Y=k) & = &
	(-1)^{k+1}
	\frac{\alpha (\alpha -1) \ldots (\alpha - k + 1)}{k!} = \nonumber\\
	& = &
	 (-1)^{k+1} { \alpha \choose k} , \qquad k = 1,2, \ldots , \label{S-mass-function}
\end{eqnarray}
and its cumulative distribution functions is
\begin{eqnarray}
	F_k = \sum_{j=1}^{k} p_j 
	& = & 
	1 - (-1)^{k} {\alpha -1  \choose k } = 
	1- \frac{1}{k \, B (k, 1-\alpha)}  =  \nonumber\\
	& = & 
	1 - \frac{\Gamma (k - \alpha + 1)}{k \, \Gamma (k) \Gamma (1-\alpha)},
	\quad
	k = 1,2, \ldots , \label{S-cumulative-function}
 \end{eqnarray}
 where $B(x,y) = \Gamma(x) \, \Gamma (y) /\Gamma (x+y)$, $x>0$, $y>0$, 
 is Euler's beta function.
 
 The probability generating function is 
 $$
	G_Y (s) = \mathbb{E} \, s^{Y} =
 	\sum_{k=1}^{\infty}
	s^{k}
	(-1)^{k+1} {\alpha \choose k} =
	1 - (1-s)^\alpha, 
	\quad 
	|s| < 1.
 $$

Note that, as $k \longrightarrow \infty$,
$$
	p_k(\alpha) \sim 
	\frac{\alpha}
		{\Gamma (1-\alpha)} \,
	\frac{1}{k^{1+\alpha}},
$$
and hence 
$$
	\mathbb{E} Y = \infty, 
$$
where $\mathbb{E}$ denotes the mathematical expectation.

The Sibuya random variable is an integer-valued random variable 
with the probability mass function given by (\ref{S-mass-function})
(see, e.g. \cite{Sibuya-1979,Sibuya-Shimitzu-1981,Devroye-1993,Pillai-Jayakumar-1995}).
Sibuya distribution can be simulated using at least three methods.


\subsection{Method 1}\label{sec:method-1}

This method has been presented and used in \cite{LP-MCFD-1}.  
We introduce $F_j = \sum_{k=1}^{j} p_k$, 
see (\ref{S-cumulative-function}), 
where $p_k$ are defined by (\ref{S-mass-function}). 
Then 
$0=F_0 < F_1 < \ldots < F_j < \ldots$, and $p_j = F_j - F_{j-1}$. 

If $U$ is a random variable uniformly distributed in $(0, 1)$, then 
$$
	\mathbb{P} (F_{j-1} < U < F_j) = p_j,
$$
and hence to generate the Sibuya random variable $Y \in \{ 1, 2, \ldots \}$, we set 
$$
Y = k,   \quad  \mbox{if~  } F_{k-1} < U < F_k.
$$.


\subsection{Method 2}\label{sec:method-2}

Using the definition of the Sibuya distribution we can generate 
the sequence of independent uniformly distributed random variables 
$U_1$, $U_2$, \ldots , $U_k$, \ldots, 
and for a fixed $\alpha \in (0,1)$ we have the following algorithm.

 \textit{Step 1.} Put $Y=1$, if $U_1<\alpha$; otherwise, go to Step 2. 
 
 \textit{Step 2.}  Put $Y=2$, if $U_1 \geq \alpha$, but $U_2 < \frac{\alpha}{2}$;  otherwise, go to Step 3. 
 
 \textit{Step 3.} Put $Y=3$, if $U_1 \geq \alpha$, $U_2 \geq \frac{\alpha}{2}$, but $U_3<\frac{\alpha}{3}$;  otherwise, go to... 
 
and so forth.  

In other words, we put $Y=k$ ($k \in \{ 1, 2, \ldots \}$) if 
$U_1 \geq \alpha$, $U_2 \geq \frac{\alpha}{2}$, \ldots, $U_{k-1} \geq \frac{\alpha}{k-1}$, but $U_k< \frac{\alpha}{k}$.


\subsection{Method 3}\label{sec:method-3}

We consider Sibuya distribution, that is a random variable $Y \in {1, 2, 3, \ldots}$, 
such that 
\begin{displaymath}
 	p_k = \mathbb{P} (Y=k) = 
	(-1)^{k+1}
	\frac{\Gamma (\alpha + 1)}{k! \, \Gamma (\alpha - k + 1)}, 
	\quad
	k = 1, 2, 3, \ldots
\end{displaymath}

Note that its probability generating function is 
\begin{equation} \label{Sibuya-pgf}
	G_Y (s) = \mathbb{E} s^{Y} 
	              = \sum_{k=1}^{\infty} s^k p_k
		      = 1 - (1-s)^\alpha, 
	\quad 
	0 < s < 1.
\end{equation}

Let us also consider a random variable $X \in {1, 2, 3, \ldots}$ with geometric distribution
\begin{equation}\label{geom-distribution}
	\mathbb{P} (X = k) = q^{k-1} p, 
	\quad
	k = 1, 2, 3, \ldots
\end{equation}
\noindent
and probability generating function
\begin{equation} \label{geom-pgf}
	G_X (s) = \mathbb{E} s^X  = 
			p s \sum_{k=1}^{\infty} (qs)^{k-1} 
			= 
			\frac{ps}{1- qs}, 
			\quad
			0 < s < 1/q. 
\end{equation}

Also we consider the beta distribution $B (\alpha, 1-\alpha)$ with density
\begin{equation} \label{beta-pdf}
	f_{\alpha} (x) = 
	\left\{
		\begin{array}{ll}
			\frac{x^{\alpha-1}(1-x)^{-\alpha}}
			       { \Gamma(\alpha) \, \Gamma(1-\alpha) } , 
			& \quad x \in {0, 1} \\
			0, & \quad x  \notin (0, 1)
		\end{array}
	\right.
\end{equation}

Note that for any fixed $\alpha \in (0, 1)$ we have 
$\lim_{x \rightarrow 0+} f_{\alpha} (x) = \infty$, 
and 
$\lim_{x \rightarrow 1-} f_{\alpha} (x) = \infty$. 
The particular case when $\alpha = \frac{1}{2}$ is known as $\arcsin$ distribution with the probability generating function 
\begin{equation}
	f_{\alpha} (x) = 
	\left\{
		\begin{array}{ll}
			\frac{1}
			       { \pi \sqrt{x (1-x)} } , 
			& \quad x \in (0, 1) \\
			0, & \quad x  \notin (0, 1)
		\end{array}
	\right.
\end{equation}

Its cdf is of the form 
\begin{equation}
	F(x) = \frac{2}{\pi} \arcsin (\sqrt{x})
               = 
	\left\{
		\begin{array}{ll}
			0, & \quad x < 0 \\
			\frac{\arcsin (2x-1)}
			       { \pi } , 
			& \quad x \in (0, 1) \\
			1, & \quad x > 1
		\end{array}
	\right.
\end{equation}

Sibuya distribution can be interpreted as randomly mixed geometric distribution.
This has been mentioned for Sibuya distribution for the case of support 
on \{ 0, 1, 2, \ldots \} by Sibuya himself \cite{Sibuya-1979} 
and later by Kozubowski and Podgorski \cite{KP-2018}. 
Since in the above papers the full proof is not presented, 
we present the proof here.

\begin{thm} \label{Sibuya-theorem}
Let $B$ has a beta distribution $B(\alpha, 1-\alpha)$.
Further, conditionally $B = p$,  
which  has a geometric distribution (\ref{geom-distribution})
$$
\mathbb{P} (Y = k | B = p) = p (1-p)^{k-1}, 
\quad
k = 1, 2, 3, \ldots
$$
Then, unconditionally, $Y$ has a Sibuya distribution 
with parameter $\alpha \in (0, 1)$. 
\end{thm}

\begin{pf}
Using conditioning argument, we have from 
(\ref{Sibuya-pgf}), (\ref{geom-distribution}), (\ref{geom-pgf}), and (\ref{beta-pdf}):
\begin{eqnarray*}
\mathbb{E} s^Y 
& = & 
\int_{0}^{1} \mathbb{E} (s^Y | B = x) p(x) dx = \\
& = & 
\int_{0}^{1} 
	\sum_{k=1}^{\infty}
	\Bigl(  (1-x)s   \Bigr)^{k-1}
	\frac{x^{\alpha -1} \, (1-x)^{-\alpha}}
	       {\Gamma (\alpha) \, \Gamma(1-\alpha)} 
	dx  = \\
& = & 
\frac{1}{\Gamma (\alpha) \, \Gamma(1-\alpha)}
\sum_{k=1}^{\infty}
s^k
\int_{0}^{1} 
x^{\alpha + 1 -1} \, (1-x)^{k-\alpha-1}
dx = \\
& = & 
\frac{1}{\Gamma (\alpha) \, \Gamma(1-\alpha)}
\sum_{k=1}^{\infty}
s^k
\frac{ \Gamma (\alpha + 1) \, \Gamma (k-\alpha)}
       {k! } = \\
& = & 
\sum_{k=1}^{\infty}
s^k 
(1-\alpha) 
\left (1 - \frac{\alpha}{2} \right) 
(1 - \frac{\alpha}{3}) 
\ldots
(1 - \frac{\alpha}{k-1}) 
\frac{\alpha}{k} = \\
& = & 
1 - (1 - s)^\alpha, 
\quad 
0 < s < 1.
\end{eqnarray*}
\end{pf}

From the Theorem \ref{Sibuya-theorem}, we arrive at the following
method of simulation of Sibuya distribution. 

\textit{Step 1.}  
For a given $\alpha \in (0, 1)$, simulate the random variable $B$
with the beta distribution $B (\alpha, 1-\alpha)$, 
denoting the resulting valueof $B$ by $p \in (0, 1)$. 

\textit{Step 2.} 
 For $p \in (0, 1)$  arising from Step 1, 
simulate geometric random variable $X$ 
with probability mass function (\ref{geom-distribution}) 
and support $\{1, 2, 3, \ldots  \}$, 
denoting the resulting integer by $k \in \{1, 2, 3, \ldots  \}$. 

\textit{Step 3.}  
From the Theorem \ref{Sibuya-theorem} we obtain that Sibuya random variable 
$Y = k$.

Note that for the particular case of $\alpha = \frac{1}{2}$ 
the simulation formula is more elegant, since in this case we have 
$$
	B = \sin^2 (\frac{\pi}{2} U),
$$
where $U$ has uniform distribution on $(0, 1)$.

\subsection{Numerical simulation of Sibuya distribution} \label{sec:num-sim-Sibuya}

The absence of tools for the simulation of the Sibuya distribution motivated us
to create the MATLAB implementation~\cite{Podlubny-Sibuya-toolbox} of all three methods, 
described in sections \ref{sec:method-1}, \ref{sec:method-2}, and \ref{sec:method-3}. 

\begin{figure}[p]

\includegraphics[width=0.3\textheight,angle=-90, bb=100 10 1120 410, clip]{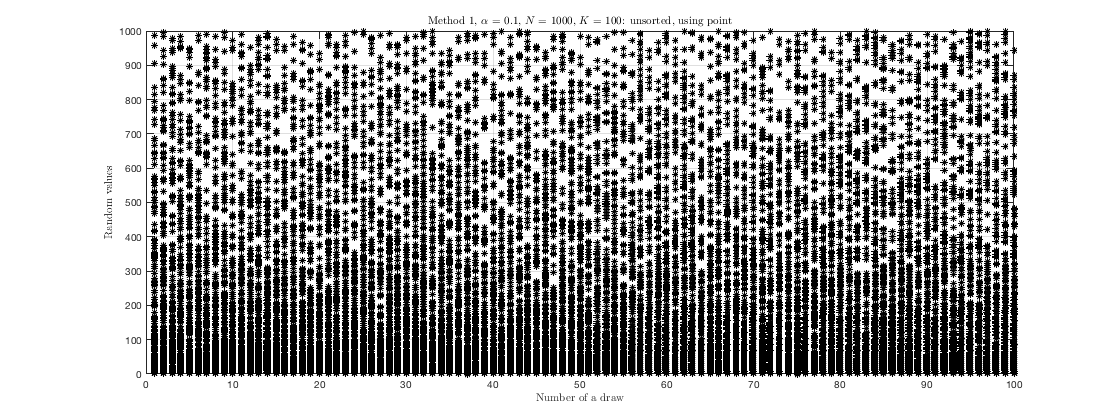}
\includegraphics[width=0.3\textheight,angle=-90, bb=100 10 1120 410, clip]{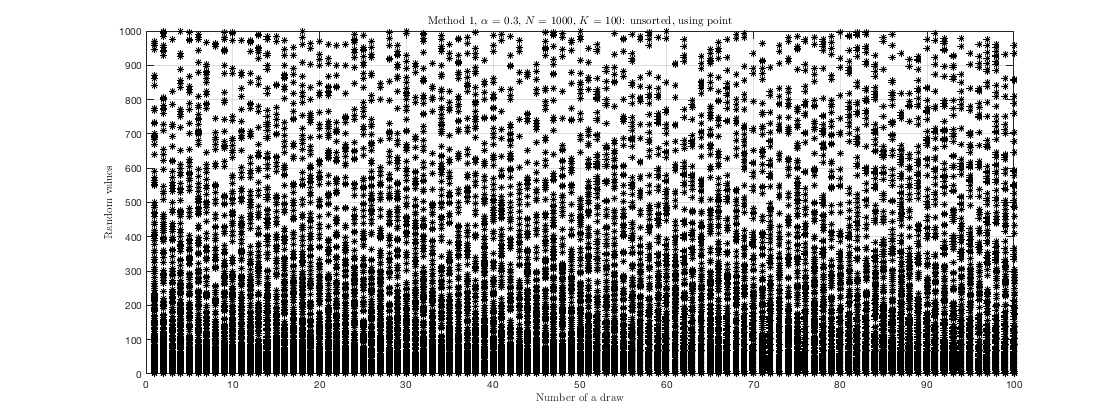}
\includegraphics[width=0.3\textheight,angle=-90, bb=100 10 1120 410, clip]{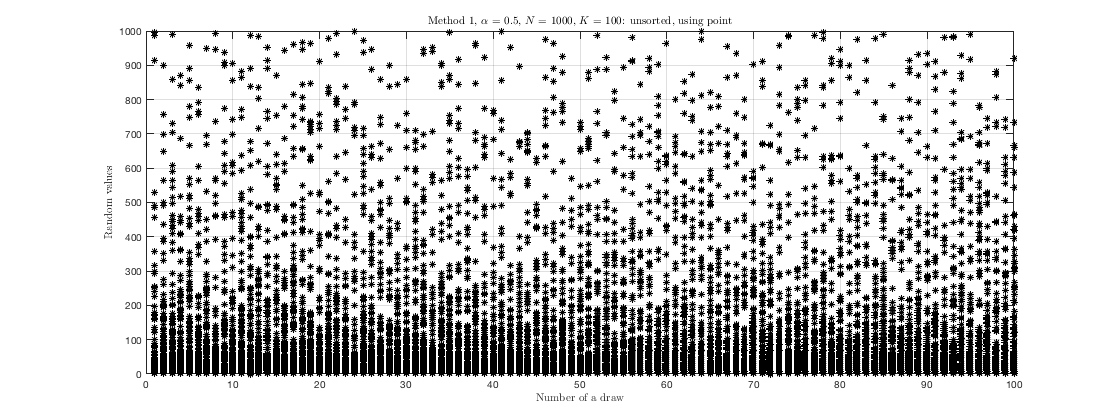}
\includegraphics[width=0.3\textheight,angle=-90, bb=100 10 1120 410, clip]{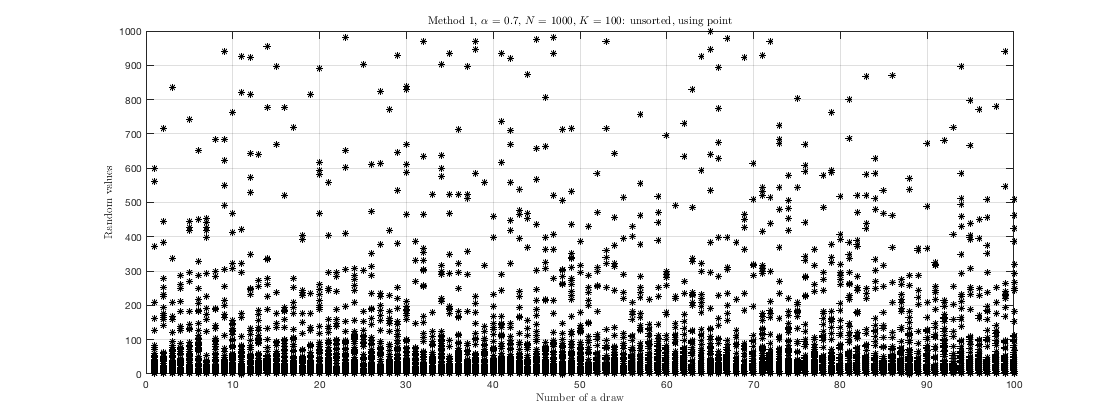}
\includegraphics[width=0.3\textheight,angle=-90, bb=100 10 1120 410, clip]{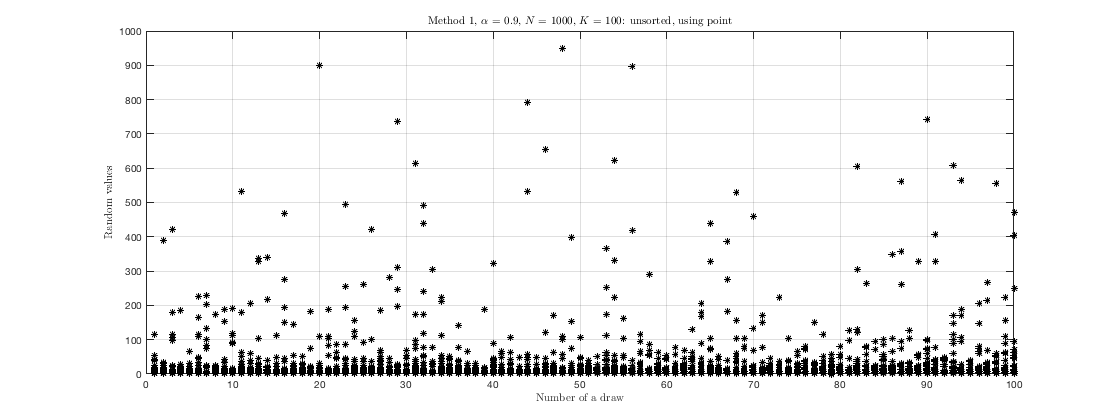}

\includegraphics[width=0.3\textheight,angle=-90, bb=100 10 1120 410, clip]{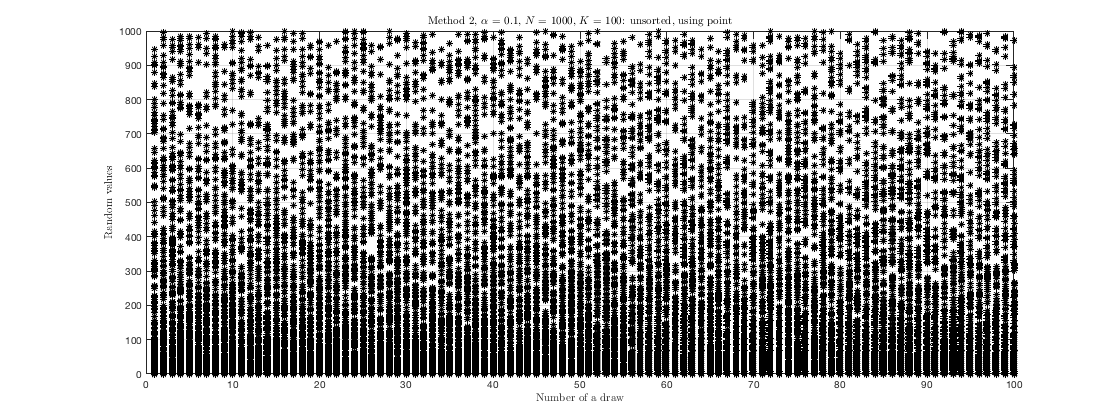}
\includegraphics[width=0.3\textheight,angle=-90, bb=100 10 1120 410, clip]{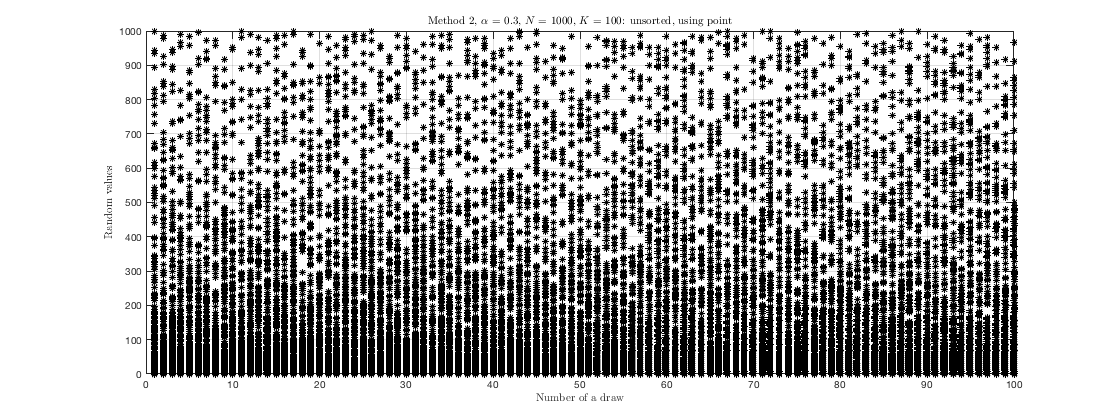}
\includegraphics[width=0.3\textheight,angle=-90, bb=100 10 1120 410, clip]{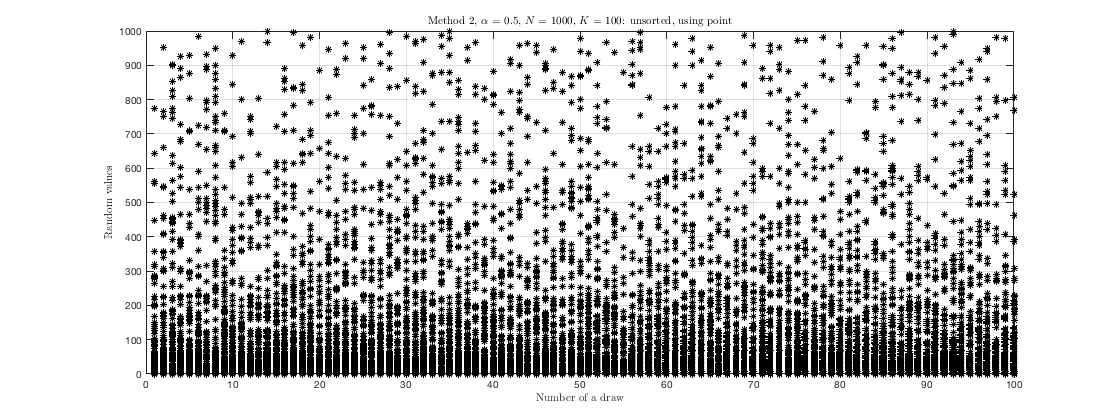}
\includegraphics[width=0.3\textheight,angle=-90, bb=100 10 1120 410, clip]{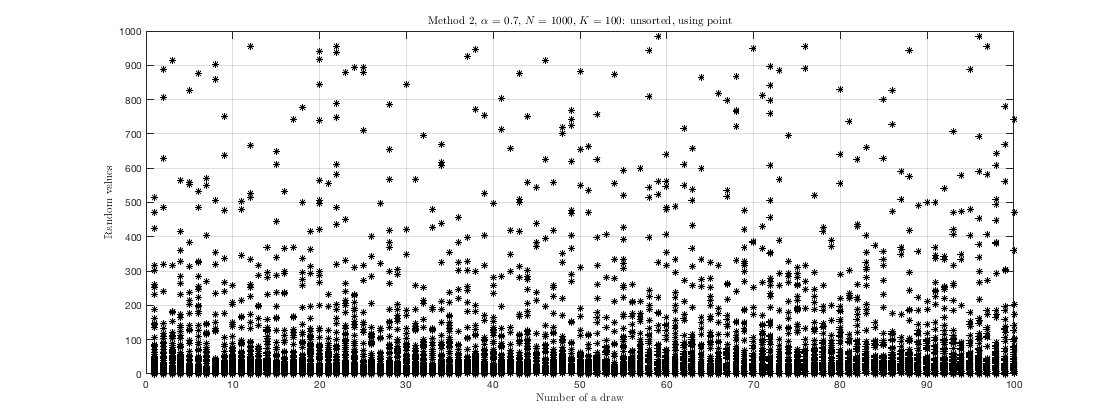}
\includegraphics[width=0.3\textheight,angle=-90, bb=100 10 1120 410, clip]{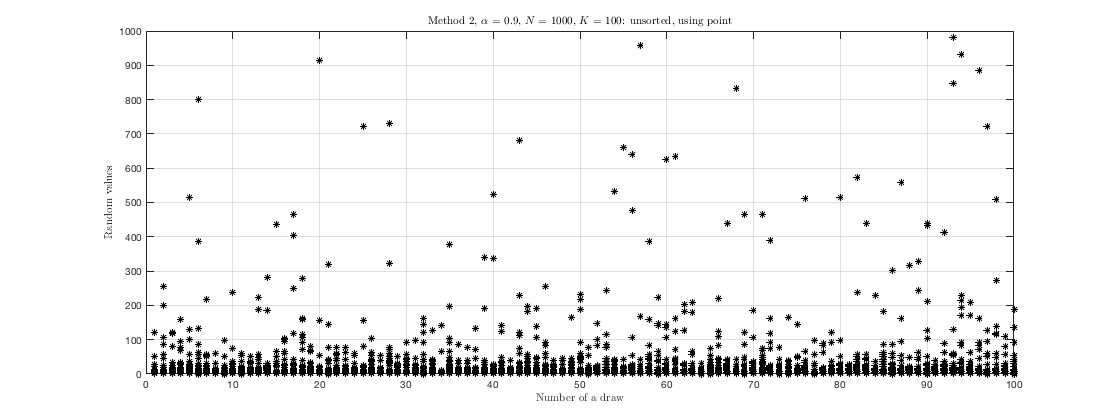}

\includegraphics[width=0.3\textheight,angle=-90, bb=100 10 1120 410, clip]{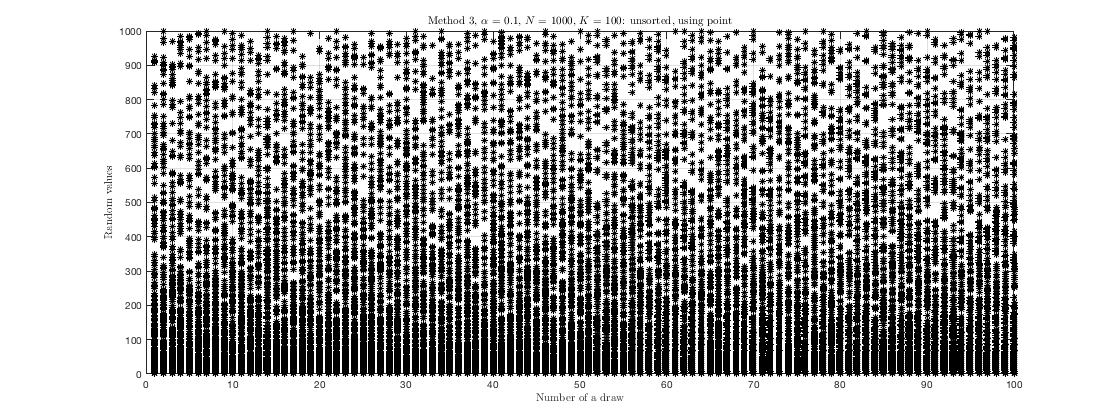}
\includegraphics[width=0.3\textheight,angle=-90, bb=100 10 1120 410, clip]{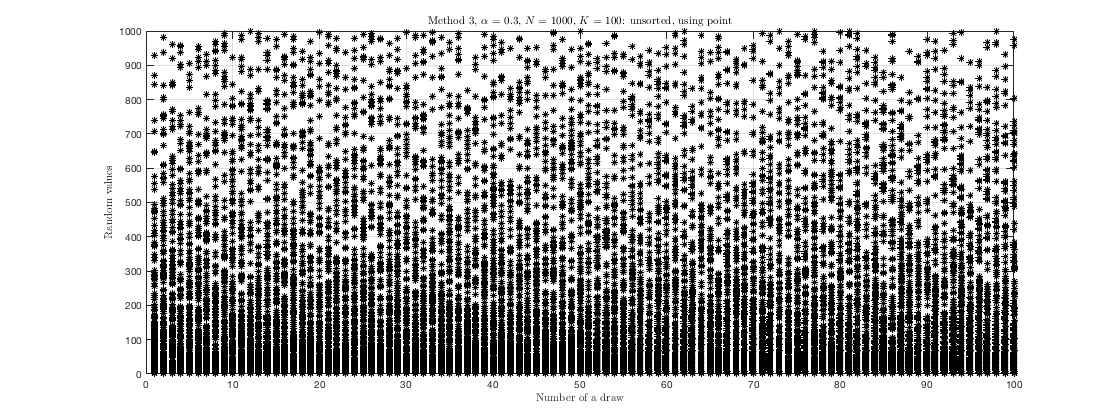}
\includegraphics[width=0.3\textheight,angle=-90, bb=100 10 1120 410, clip]{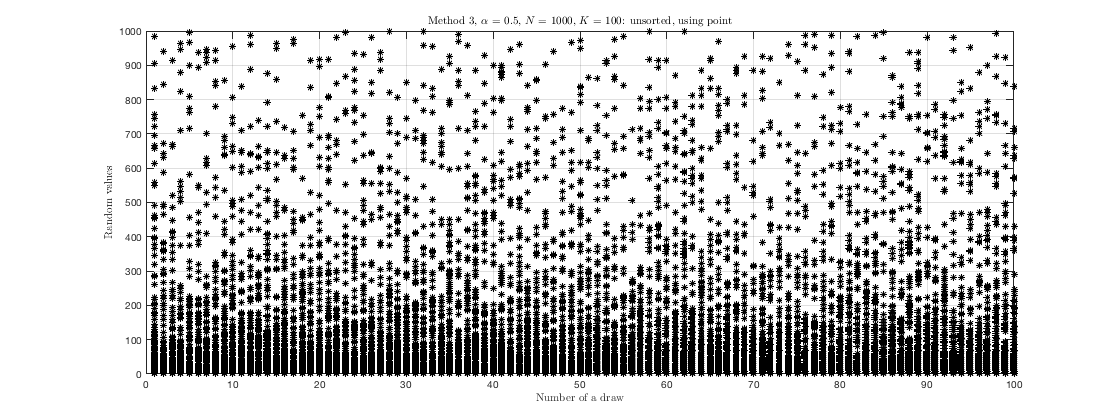}
\includegraphics[width=0.3\textheight,angle=-90, bb=100 10 1120 410, clip]{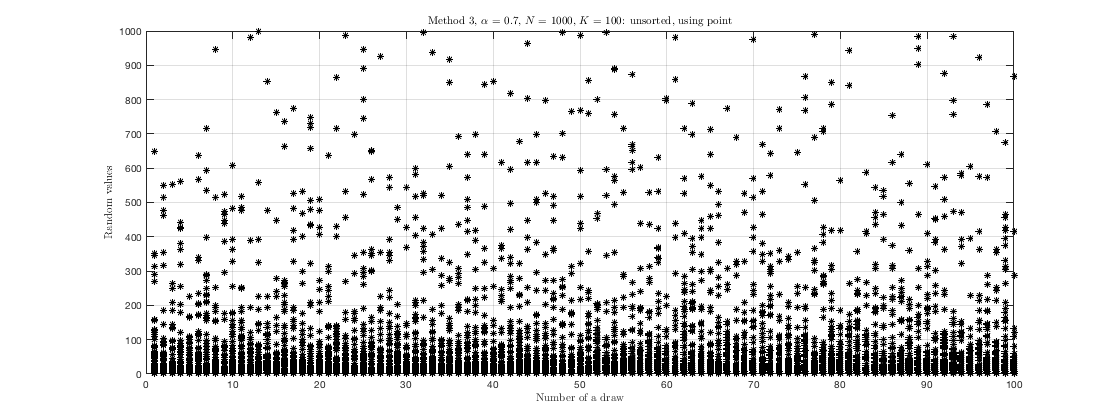}
\includegraphics[width=0.3\textheight,angle=-90, bb=100 10 1120 410, clip]{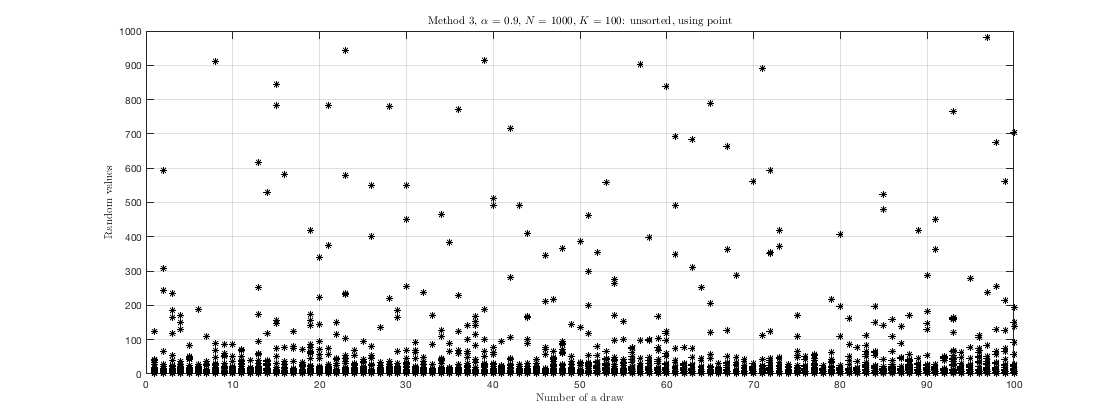}

\caption{Simulation of the Sibuya distribution by methods 1 (top), 2 (middle), and 3 (bottom) 
for (left to right)  $\alpha=0.1$, $\alpha=0.3$, $\alpha=0.5$, $\alpha=0.7$,  $\alpha=0.9$,
with 100 draws of 1000 numbers each.}
\label{fig:sims}
\end{figure}

The examples of the simulations, shown in Fig.~\ref{fig:sims}, 
provide the visualization of the Sibuya distribution for  $\alpha=0.1$, $\alpha=0.3$, 
$\alpha=0.5$, $\alpha=0.7$,  and  $\alpha=0.9$. 
All three methods deliver the same character of the distribution of random values, as seen in  Fig.~\ref{fig:sims}. 
Each image there contains one hundred of horizontal sets of 1000 random points (one hundred draws).  
The density of the random points is higher near the left end (smaller values), 
and it decreases towards the right end (higher values).  
For higher values of $\alpha$ the random points concentrate or gravitate, so to say,  toward the left end. 

\begin{figure}[htbp]
\begin{center}

\includegraphics[width=0.3\textwidth, bb=40 10 520 410, clip]{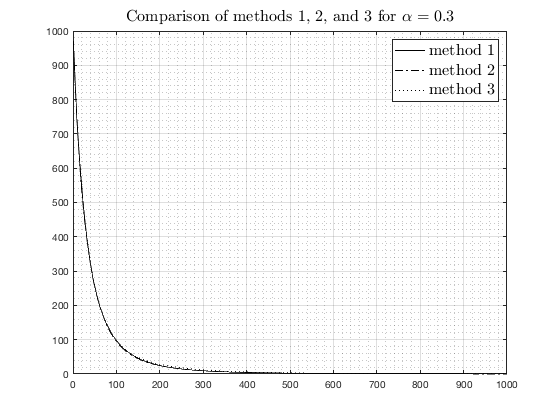}
\hfill
\includegraphics[width=0.3\textwidth, bb=40 10 520 410, clip]{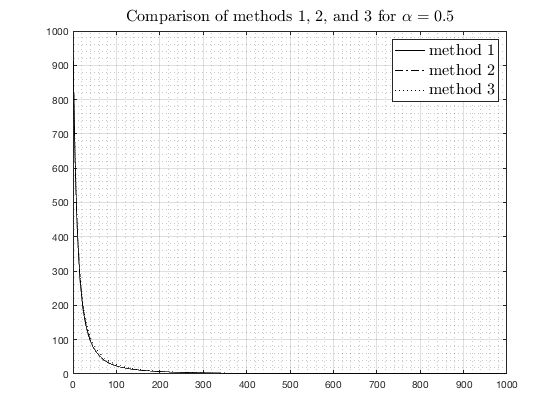}
\hfill
\includegraphics[width=0.3\textwidth, bb=40 10 520 410, clip]{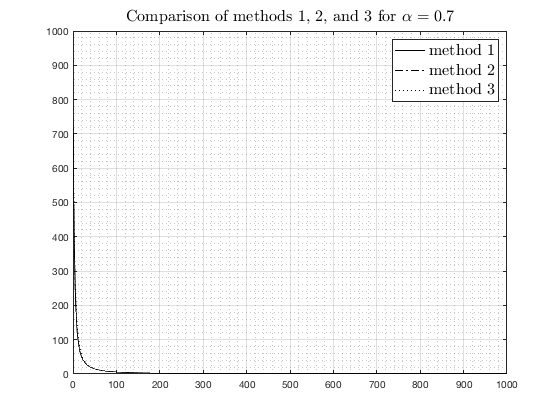}

\caption{Comparison of the methods 1, 2, and 3 using the sorted averages 
of the the drawn random numbers for $\alpha=0.3$ (left), $\alpha=0.5$ (middle), $\alpha=0.7$ (right).}
\label{fig:compare}
\end{center}
\end{figure}

The three methods are compared using the sorted averages 
of the the drawn random numbers. To do this, the random numbers of each draw are first sorted, 
and then average of all draws is computed. The results for $\alpha=0.3$, $\alpha=0.5$, $\alpha=0.7$ 
are shown in Fig.~\ref{fig:compare}. 

The presented three methods of different nature deliver practically identical results,  
which shows that any of them can be used for simulating the Sibuya distribution.

\section{Sibuya-like distributions and non-local operators} \label{sec:4}

\subsection{The general Monte Carlo procedure}

We consider again the non-local operator (\ref{GL-definition}), 
which is known as the Gr\"unwald--Letnikov fractional derivative
and could be used (under certain conditions) also for numerical approximation 
of the Caputo  and the Riemann--Liouville fractional derivatives
(see, for example, \cite{Podlubny-FDE,LP-MCFD-1,LP-MCFD-2} ). 

We propose a general Monte Carlo scheme for numerical approximations of the operator
\begin{equation}\label{GL-operator}
 \Delta_h^\alpha f(t) = \sum_{k=0}^{\infty} w_k \, f(t - kh), 
 \quad 
 t > 0, \, \alpha > 0,
\end{equation}
where
\begin{equation}\label{GL-weights}
	w_k = (-1)^k
	\frac{ \Gamma (\alpha + 1)}
	       {k! \, \Gamma (\alpha - k + 1)}, 
	 \quad
	 k = 0, 1, 2, \ldots
\end{equation}

Noting that $w_0 = 1$, and having for $|s| \leq 1$  the following 
binomial series expansion:
$$
	(1 - s)^\alpha = \sum_{k=0}^\infty w_k \, s^k,
$$ 
we by putting $s = 1$ have 
$$
   \sum_{k=0}^\infty w_k  = 1 + \sum_{k=1}^\infty w_k = 0, 
$$
or
\begin{equation}
	\sum_{k=1}^\infty w_k = -1, 
	\quad
	\sum_{k=1}^\infty (-w_k) = 1,
	\quad
	\sum_{k=1}^\infty \gamma_k = 1;
\end{equation}
but the numbers $\gamma_k = -w_k$ could not be always positive (see \cite{LP-MCFD-2}).

In order to overcome these difficulties we introduce
$$
	I_+ = \{ k \geq 1 :   w_k > 0 \},
$$
where we assume that $|I_+| \geq 2$, or $I_+ = \varnothing$, and 
$$
	W_+ = \sum_{k \in I_+} w_k,
$$
then on the discrete set $I_+$,  which contains at least two points, 
we can introduce the discrete distribution 
$$
   P  = \left\{  \frac{w_k}{W_+} , \,\, k \in I_+ \right\} 
$$
where obviously
$$
	p_k > 0, 
	\quad 
	\sum_{ k \in I_+}  p_k = 1.
$$

Note that in some cases all $w_k<0$, $k \in \{  1, 2, \ldots, n, \ldots \} \setminus \{ k^*\}$, 
that is there exists a single integer $k^*$: $w_{k^*}>0$. 
Thus, if $w_{k^*}=0$, then the cardinality of the set $I_{+}$ 
is always greater than or equal two: $| I_{+} | \geq 2$.

Similarly, we consider 
$$
	I_{-} = \{  k\geq 1: w_k<0  \}, 
	\mbox{ where } 
	| I_{-} | \geq 2, 
	\mbox{ or } 
	I_{-} = \emptyset,
$$
and we introduce 
$$
	W_{-} = \sum_{k\in I_{-}}^{} w_k.
$$

Then on the discrete set $I_{-}$  containing at least two points
we can introduce the discrete distribution
$$
	Q = \left\{ 
			q_k = \frac{w_k}{W_{-}}, 
			\
			k \in I_{-}
		\right\}
$$
with
$$
	q_k > 0, 
	\quad
	\sum_{k \in I_{-}}^{} q_k = 1. 
$$

Again, it is possible that all 
$w_k > 0$, $k \in \{  1, 2, \ldots, n, \ldots \} \setminus \{ k_*\}$,
where  $k_*$ is a single integer such that: $w_{k_*}< 0$. 
Thus, if $w_{k_*} = 0$, then $| I_{-} | \geq 2$.

Now we can represent the operator (\ref{GL-operator}) as follows:
\begin{eqnarray}
 \Delta_h^\alpha f(t)  & = & \sum_{k=0}^{\infty} w_k \, f(t - kh) = \nonumber \\
 & = & 
 f(t) + \sum_{k=1}^{\infty} w_k \, f(t - kh) =   \nonumber \\
 & = & 
 f(t) + W_{+} \sum_{k \in I_{+}} p_k f(t - kh) + w_{k_*} f(t - k_* h) +  \nonumber \\
 & ~ & 
 + \ W_{-} \sum_{k \in I_{-}} q_k f(t - kh) + w_{k^*} f(t - k^* h) = \nonumber \\
& = & 
f(t) + w_{k_*} f(t - k_* h)  + w_{k^*} f(t - k^* h) +  \nonumber \\
 & ~ & 
+ \ W_{+} \mathbb{E}_P f(t - Y^{+} h) + W_{-} \mathbb{E}_Q f(t - Y^{-} h),
	\label{NL-4.4}
\end{eqnarray}
where $w_{k_*} < 0$, $W_{-} < 0$, 
and assuming that the above expectations exist, 
$Y^{+}$ is discrete random variable with the support $I_{+}$
and probability mass function $\{ p_k, \, k \in I_{+}\}$, 
while $Y^{-}$ is discrete random variable with the support $I_{-}$
and probability mass function $\{ q_k, \, k \in I_{-}$\}. 
We always assume that random variables $Y^{+}$ and $Y^{-}$ are independent.

Let $Y_1^{+}$, $Y_2^{+}$, \ldots , $Y_n^{+}$, \ldots are independent copies 
of random variable $Y^{+}$,
and $Y_1^{-}$, $Y_2^{-}$, \ldots , $Y_n^{-}$, \ldots are independent copies 
of random variable $Y^{-}$. 
If~$\mathbb{E}_{P} f(t - Y^{+} h) < \infty$, $\mathbb{E}_{Q} f(t - Y^{-} h) < \infty$,
then by the strong law of the large number with probability one 
\begin{equation} \label{eq:4.5}
 \frac{1}{N}
 \sum_{n=1}^{N}  f(t - Y_n^{+} h)
 \longrightarrow 
 \mathbb{E}_{P} f(t - Y^{+} h),
\end{equation}
\begin{equation} \label{eq:4.6}
 \frac{1}{N}
 \sum_{n=1}^{N}  f(t - Y_n^{-} h)
 \longrightarrow 
 \mathbb{E}_{Q} f(t - Y^{-} h),
\end{equation}
as $N \rightarrow \infty$.

Consequently, with probability one 
\begin{eqnarray}
\left(  A_{N,h}^\alpha f   \right) (t) & = &
   \frac{1}{h^\alpha}
	\Bigl[
	  f(t) + w_{k_*} f(t - k_* h) + w_{k^*} f(t - k^* h)  
	  \Bigr.  \nonumber\\
	 & ~ &
	 + \ W_{+}  \frac{1}{N} \sum_{n=1}^{N} f(t - Y_n^{+} h) \nonumber\\
	  & ~ &
	 \Bigl. 
	  + \ W_{-}  \frac{1}{N} \sum_{n=1}^{N} f(t - Y_n^{-} h) 
	\Bigr]
	\longrightarrow
	A_h^\alpha f(t), 
	\quad N \rightarrow \infty, \label{eq:4.7}
\end{eqnarray}
for any $\alpha > 0$, $h>0$. 

The above procedure is the Monte Carlo scheme for (possible) signed measures
$W = \{ w_k, k = 1, 2, \ldots \}$ (see \cite{LP-MCFD-2}). 

We present illustrative examples for the Monte Carlo scheme (\ref{eq:4.7})
for $0 < \alpha < 1$, $1 < \alpha < 2$, $2 < \alpha < 3$, $3 < \alpha < 4$, 
and $4 < \alpha < 5$. 
In particular, the case of $0 < \alpha < 1$ is important for ultraslow diffusion, 
$1 < \alpha < 2$ is appears in fractional-order diffusion-wave problems, 
$2 < \alpha < 3$ is related to fractional-order modeling of viscoelastic rods and plates, 
while other cases are presented for completeness.

\subsection{The case of $0 < \alpha < 1$}

\begin{figure}
\begin{center}
	\includegraphics[scale=0.2]{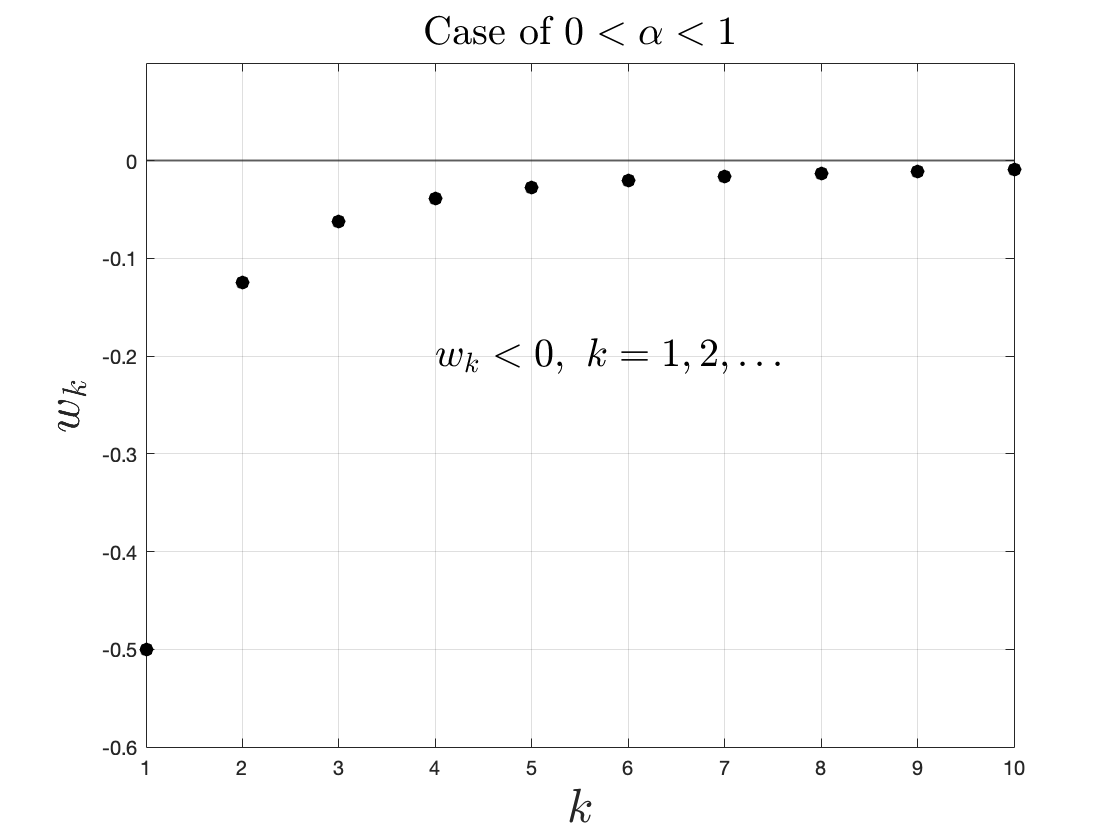}\\	
\caption{The values of $w_k$, $k = 1, 2, \ldots$, in case $0<\alpha<1$, with $\alpha = 0.5$}
\label{fig:case01}
\end{center}
\end{figure}

In this case all $w_k < 0$, $k = 1, 2, \ldots$ (see Fig.~\ref{fig:case01}), 
thus
$I_{+} = \varnothing$, $w_{k^*} = 0$, thus $Y^{+} \equiv 0$, 
while
$I_{-} = \{  k = 1, 2, \ldots: w_k < 0 \} = \{ 1, 2, 3, \ldots, n, \ldots \}$, 
and $w_{k_*} = 0$. 
The discrete random variable $Y^{-}$ has the support $\{ 1, 2, \ldots, n, \ldots \}$,
and distribution $Q$ is Sibuya distribution (see Section \ref{sec:3})
with state probabilities
\begin{equation} \label{eq:4.8}
  q_k = (-1)^{k+1}  
            \frac{\Gamma (\alpha + 1)}
                   {k! \, \Gamma (\alpha - k + 1)},
                   \quad
                   k = 1, 2, \ldots,
\end{equation}
since 
$$
   W_{-} = \sum_{k \in I_{-}} w_k = 1. 
$$

The Monte Carlo scheme (\ref{eq:4.7}) for $0 < \alpha < 1$ takes the form:
\begin{displaymath}
	\frac{1}{h^\alpha}
	\left[
		f(t) - \frac{1}{N} 
		\sum_{n=1}^{N} f(t - Y_n^{-} h)
	\right]
	\longrightarrow
	A_h^\alpha f(t)
\end{displaymath}
with probability one as $N \rightarrow \infty$, 
where $Y_1^{-}$, $Y_2^{-}$, \ldots , $Y_n^{-}$, \ldots are independent copies 
of Sibuya distribution (\ref{S-mass-function}).

\subsection{The case of $1 < \alpha < 2$}

\begin{figure}
\begin{center}
	\includegraphics[scale=0.2]{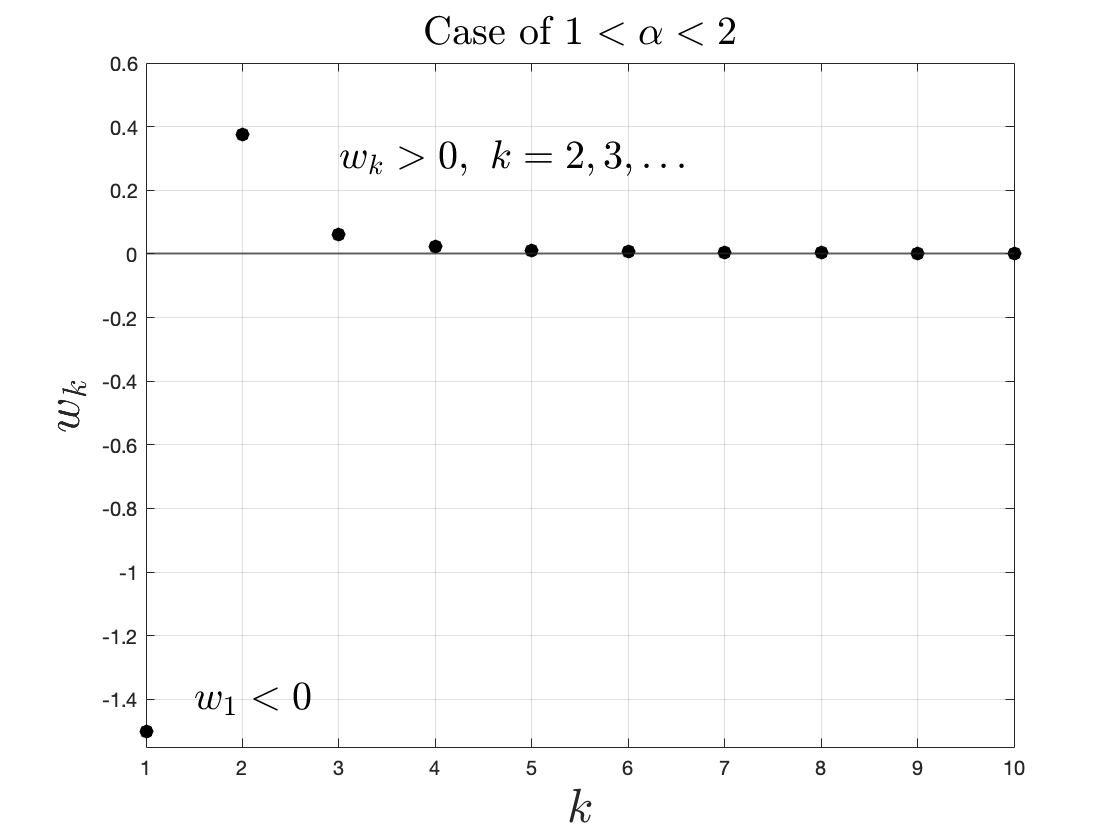}\\	
\caption{The values of $w_k$, $k = 1, 2, \ldots$, in case $1<\alpha<2$, with $\alpha = 1.5$}
\label{fig:case12}
\end{center}
\end{figure}

In this case 
$w_1 = - \alpha < 0$, 
$w_k > 0$, $k = 2, 3, \ldots$ (see Fig.~\ref{fig:case12}),
thus
$I_{+} =  \{  k>0: w_k > 0  \} = \{ 2, 3, \ldots \} $, 
while 
$k_* = 1$, $w_1 = w_{k_*} =  -\alpha$, and 
$-\alpha + \sum\limits_{k=2}^{\infty} w_k = -1$, hence
\begin{displaymath}
	W_{+} =  \sum_{k=2}^{\infty} w_k = \alpha - 1 > 0. 
\end{displaymath}

Obviously, $Y^{-} \equiv 0$, since $I_{-} = \varnothing$, $w_{k^{*}} = 0$. 

Thus, if $Y^{+}$ is the discrete random variable with support 
$\{ 2, 3, \ldots, n, \ldots \}$, and
\begin{eqnarray}
 	p_k 
	& = &  \mathbb{P} (Y^{+} = k) = \frac{w_k}{W_{+}} = \frac{w_k}{\alpha -1} 
		\nonumber\\
	& = & 
	\frac{1}{\alpha -1} (-1)^{k} 
	\frac{\Gamma (\alpha + 1)}
	       {k! \, \Gamma (\alpha - k + 1)}, 
	 \quad
	 k = 2, 3, \ldots ,
	 		\label{eq:4.9}
\end{eqnarray}
and $Y_1^{+}$, $Y_2^{+}$, \ldots , $Y_n^{+}$, \ldots are independent copies 
of random variable $Y^{+}$ given by (\ref{eq:4.9}), 
the Monte Carlo procedure takes on the form:
\begin{displaymath}
	\frac{1}{h^\alpha}
	\Bigl[
		f(t) - \alpha f(t-h) 
		 + (\alpha -1) \frac{1}{N} \sum_{n=1}^{N} f(t - Y_n^{+} h)
	\Bigr]
	\longrightarrow
	A_h^\alpha f(t),
\end{displaymath}
with probability one, as $N \rightarrow \infty$.

\subsection{The case of  $2 < \alpha < 3$}

\begin{figure}
\begin{center}
	\includegraphics[scale=0.2]{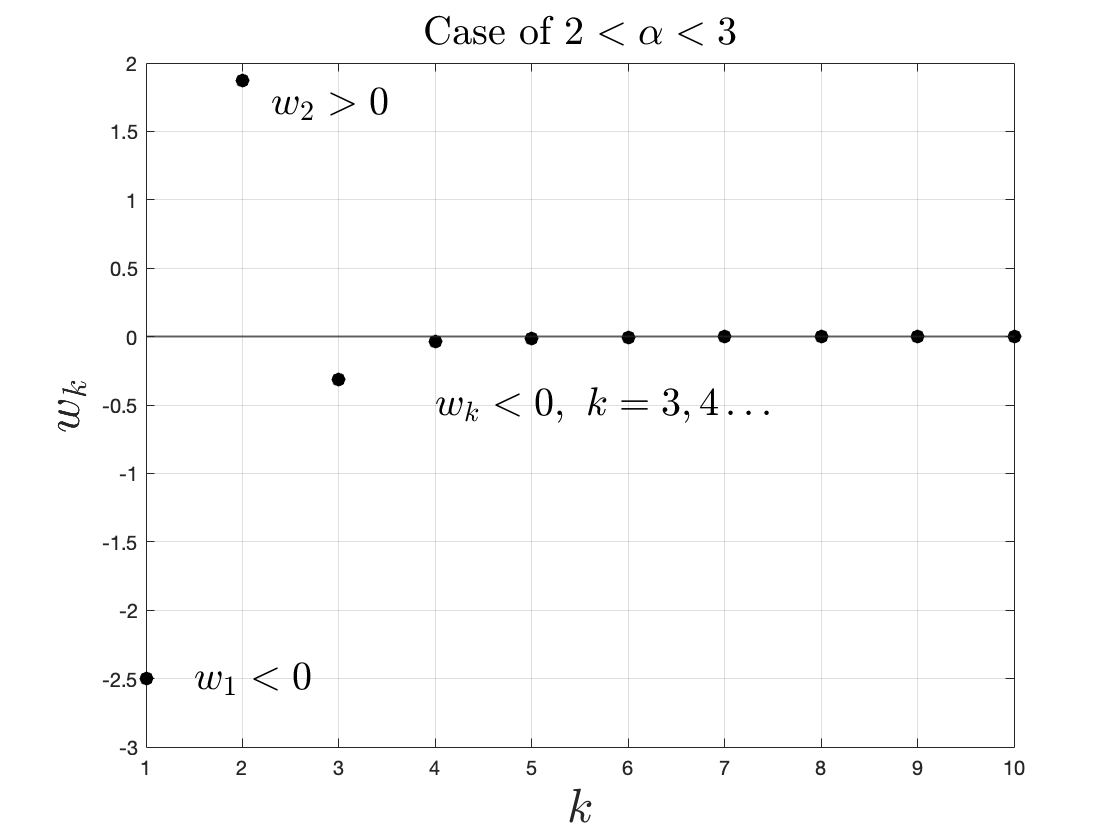}\\	
\caption{The values of $w_k$, $k = 1, 2, \ldots$, in case $2<\alpha<3$, with $\alpha = 2.5$}
\label{fig:case23}
\end{center}
\end{figure}

In this case
$w_1 = -\alpha < 0$, $w_2 = \frac{\alpha (\alpha-1)}{2} > 0$, 
$w_k < 0$, $k \geq 3$
(see Fig.~\ref{fig:case23}), 
and $Y^{+} \equiv 0$, since $I_{+} = \varnothing$, 
but $k^{*}=2$, and 
$w_{k^{*}}=\frac{\alpha (\alpha -1)}{2}$. 

Then 
$I_{-} = \{ 1, 3, 4, \ldots, n, \ldots \}$, and
\begin{displaymath}
	-\alpha + \frac{\alpha (\alpha -1)}{2} 
	+ \sum_{k=3}^{\infty} w_k = -1,
\end{displaymath}
thus 
\begin{displaymath}
	W_{-} = -\alpha + \sum_{k=3}^{\infty} w_k = -1 - \frac{\alpha (\alpha -1)}{2} < 0.
\end{displaymath}

We introduce the random variable $Y^{-}$ 
with support $I_{-} = \{ 1, 3, 4, \ldots, n, \ldots \}$ 
and state probabilities:
\begin{equation} \label{eq:4.10}
	q_1 = \frac{-\alpha}{W_{-}} = \frac{2 \alpha}{(2 + \alpha (\alpha -1))} > 0,
\end{equation}
\begin{equation} \label{eq:4.11}
	q_k = \frac{w_k}{W_{-}} 
	= 
	\frac{2 \, (-1)^{k+1} \Gamma (\alpha + 1)}
	       {(2 + \alpha (\alpha -1)) \, k! \, \Gamma (\alpha - k + 1)}, 
	\quad 
	k = 3, 4, \ldots
\end{equation}

We call the distribution $Q$ of the random variable $Y^{-}$
the \emph{sieved Sibuya-like distribution}
for $2 < \alpha < 3$.
 
 We arrive at the following Monte Carlo procedure. 
 If  $Y_1^{-}$, $Y_2^{-}$, \ldots, $Y_n^{-}$,~\ldots are independent copies of the 
 sieved Sibuya-like distribution (\ref{eq:4.10})--(\ref{eq:4.11}), 
 then with probability one as $N \rightarrow \infty$
 \begin{eqnarray*}
 	\frac{1}{h^\alpha}
 	\Bigl[
		f(t) + \frac{\alpha (\alpha -1)}{2} f(t-2h) 
		  - \left(
		  	1 + \frac{\alpha (\alpha -1)}{2}
		    \right)
		    \frac{1}{N}
		    \sum_{n=1}^{N} f(t - Y_n^{-} h)
	\Bigr]
	\rightarrow
	A_h^{\alpha} f(t)
 \end{eqnarray*}

 \subsection{The case of  $3 < \alpha < 4$}

\begin{figure}
\begin{center}
	\includegraphics[scale=0.2]{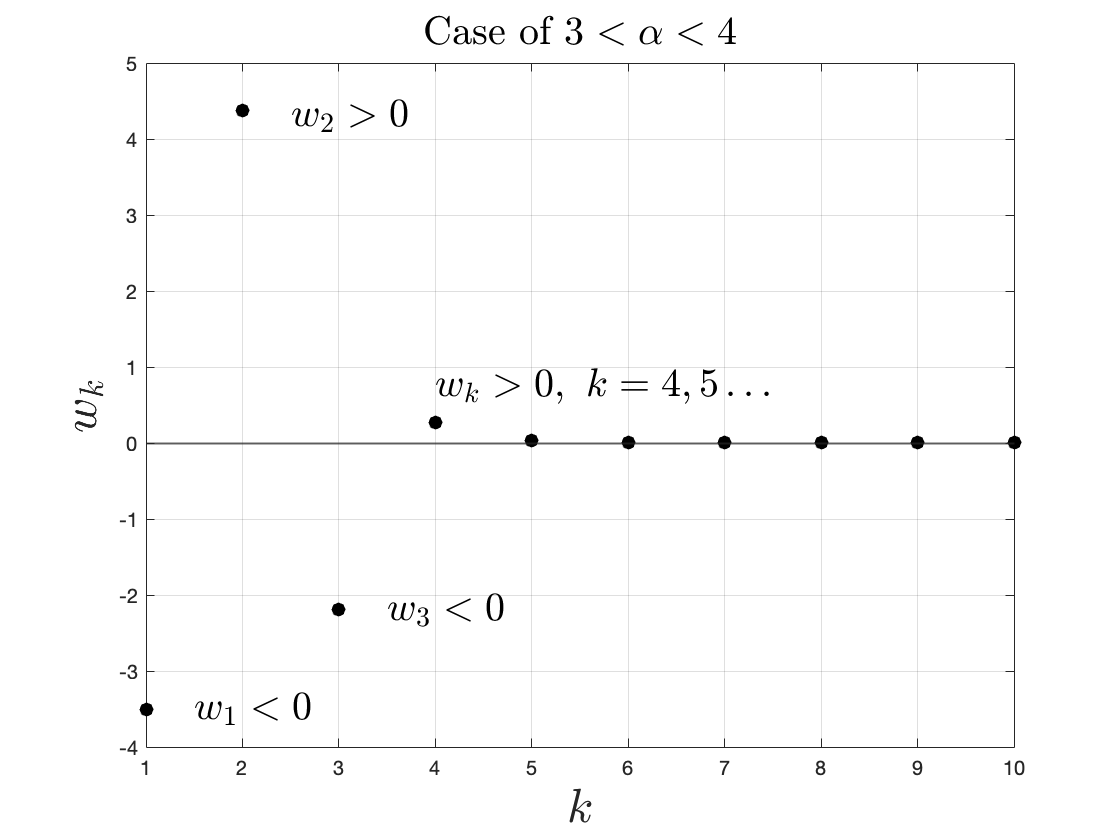}\\	
\caption{The values of $w_k$, $k = 1, 2, \ldots$, in case $3<\alpha<4$, with $\alpha = 3.5$}
\label{fig:case34}
\end{center}
\end{figure}

In this case
$w_1 = -\alpha < 0$,
$w_2 = \frac{1}{2} \frac{\Gamma (\alpha + 1)}{\Gamma (\alpha -1)} > 0$, 
$w_3 = -\frac{1}{6} \frac{\Gamma (\alpha + 1)}{\Gamma (\alpha -2)} < 0$,
$w_k > 0$, $k \geq 4$
(see Fig~\ref{fig:case34}),
and
\begin{displaymath}
	-\alpha + \frac{1}{2} \frac{\Gamma (\alpha + 1)}{\Gamma (\alpha -1)}
	-\frac{1}{6} \frac{\Gamma (\alpha + 1)}{\Gamma (\alpha -2)} 
	+ \sum_{k=4}^{\infty} w_k = -1,
\end{displaymath}
 thus
 \begin{displaymath}
 	W_{+} = w_2 + \sum_{k=4}^{\infty} w_k 
	= 
	\alpha - 1 + \frac{1}{6} \frac{\Gamma (\alpha + 1)}{\Gamma (\alpha -2)} > 0.
 \end{displaymath}

 Let $Y^{+}$ is the discrete random variable 
 with support $\{  2, 4, 5, 6, \ldots, n, \ldots\}$
 and sieved Sibuya-like distribution $P$ for $3 <\alpha < 4$,
 with state probabilities
 \begin{displaymath} 
 	p_2 = \mathbb{P} ( Y^{+} = 2 )
	= 
	\frac{\frac{1}{2}  \frac{\Gamma (\alpha + 1)}{\Gamma (\alpha -1)}}
	        {\alpha - 1 + \frac{1}{6} \frac{\Gamma (\alpha + 1)}{\Gamma (\alpha -2)} },
 \end{displaymath}
  \begin{equation}  \label{eq:4.12}
 	p_k =  \mathbb{P} ( Y^{+} = k )  
	= 
	\frac{ (-1)^k \, \Gamma (\alpha +1) \frac{1}{k! \, \Gamma (\alpha - k + 1)} }
	       {\alpha - 1 + \frac{1}{6} \frac{\Gamma (\alpha + 1)}{\Gamma (\alpha -2)} },
	 \quad 
	 k = 4, 5, \ldots
  \end{equation}
  and $Y_1^{+}$, $Y_2^{+}$, \ldots, $Y_n^{+}$, \ldots
  are independent copies of random variable $Y^{+}$.
  
  The set 
  $I_{-} = \{ k: w_k < 0 \}$ = \{ 1, 3 \},
  and $Y^{-}$ is the random variable with support $\{ 1, 3\}$
  and state probabilities 
  \begin{equation}  \label{eq:4.13}
  	q_1 = \mathbb{P} (Y^{-} = 1) 
	=
	\frac{\alpha}
	       {\alpha + \frac{1}{6} \frac{\Gamma (\alpha + 1)}{\Gamma (\alpha -2)} }, 
  \end{equation}
    \begin{equation}  \label{eq:4.14}
  	q_3 = \mathbb{P} (Y^{-} = 3) 
	=
	\frac{\frac{1}{6} \frac{\Gamma (\alpha + 1)}{\Gamma (\alpha -2)} }
	       {\alpha + \frac{1}{6} \frac{\Gamma (\alpha + 1)}{\Gamma (\alpha -2)} }, 
  \end{equation}
 since 
 $W_{-} = -\alpha - \frac{1}{6} \frac{\Gamma (\alpha + 1)}{\Gamma (\alpha -2)}  < 0$.
 
 Note that for $3 < \alpha < 4$ we have 
 $w_{k_{*}} = w_{k^{*}}=0$, 
 that is there is no single positive or negative probabilities.

 Let $Y_1^{-}$, $Y_2^{-}$, \ldots, $Y_n^{-}$, \ldots 
 are independent copies of random variable $Y^{-}$, 
 then as $N \rightarrow \infty$ with probability one we have 
 the following Monte Carlo procedure:
 \begin{eqnarray*}
 	& ~ &
 	\frac{1}{h^\alpha}
	\left[
	    f(t) 
	    + \left(
	    	\alpha -1  + \frac{1}{6} \frac{\Gamma (\alpha + 1)}{\Gamma (\alpha -2)} 
	    \right)
	    \frac{1}{N}
	    \sum_{n=1}^{N} f(t - Y^{+} h)
	  \right.
	    		\\
 	& ~ &		
	\hspace*{4em}	
	  \left.
	    - \left(
	         \alpha + \frac{1}{6} \frac{\Gamma (\alpha + 1)}{\Gamma (\alpha -2)} 
	      \right)
	      \frac{1}{N}
	       \sum_{n=1}^{N} f(t - Y^{-} h)
	\right]
	\longrightarrow 
	A_h^{\alpha} f(t).
 \end{eqnarray*}

  \subsection{The case of  $4 < \alpha < 5$}
  
\begin{figure}
\begin{center}
	\includegraphics[scale=0.2]{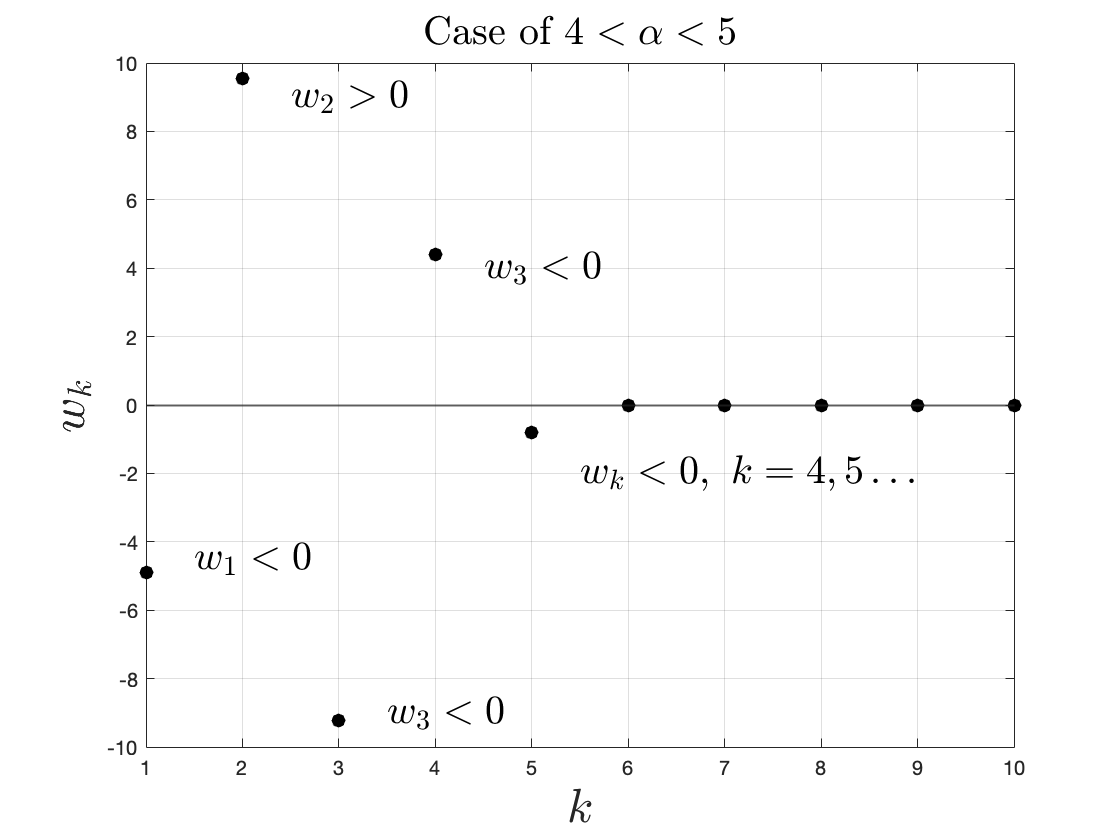}\\	
\caption{The values of $w_k$, $k = 1, 2, \ldots$, in case $4<\alpha<5$, with $\alpha = 4.5$}
\label{fig:case45}
\end{center}
\end{figure}

  In this case
  $w_1 = -\alpha < 0$,
  $w_2 = \frac{1}{2} \frac{\Gamma (\alpha + 1)}{ \Gamma (\alpha -1)} > 0$,
  $w_3 = -\frac{1}{6} \frac{\Gamma (\alpha + 1)}{\Gamma (\alpha -2)} < 0$,
  $w_4 = \frac{1}{24} \frac{\Gamma (\alpha + 1)}{\Gamma (\alpha -3)} > 0$, 
  $w_k < 0$, $k = 5, 6, \ldots$
  (see Fig.~\ref{fig:case45}).
  
 Thus,  $I_{+} = \{ 2, 4 \}$,
 $$
 	W_{+} = \frac{1}{2} \frac{\Gamma (\alpha + 1)}{ \Gamma (\alpha -1)}
	             + \frac{1}{24} \frac{\Gamma (\alpha + 1)}{\Gamma (\alpha -3)}  > 0, 
 $$
 and random variable $Y^{+}$ has the support $I_{+} = \{ 2, 4 \}$
 with probabilities
 \begin{equation} \label{eq:4.15}
 	p_2 = \mathbb{P} (Y^{+} = 2) 
	       = 
	       \frac{\frac{1}{2} \frac{\Gamma (\alpha + 1)}{ \Gamma (\alpha -1)}}
	              { \frac{1}{2} \frac{\Gamma (\alpha + 1)}{ \Gamma (\alpha -1)}
	             + \frac{1}{24} \frac{\Gamma (\alpha + 1)}{\Gamma (\alpha -3)}},
 \end{equation}
  \begin{equation} \label{eq:4.16}
 	p_4 = \mathbb{P} (Y^{+} = 4) 
	       = 
	       \frac{\frac{1}{24} \frac{\Gamma (\alpha + 1)}{\Gamma (\alpha -3)}}
	              { \frac{1}{2} \frac{\Gamma (\alpha + 1)}{ \Gamma (\alpha -1)}
	             + \frac{1}{24} \frac{\Gamma (\alpha + 1)}{\Gamma (\alpha -3)}},
 \end{equation}

The random variable $Y^{-}$ is the sieved Sibuya-like distribution
for $4 < \alpha < 5$ with its support
$$
	I_{-} = \{  1, 3, 5, 6, \ldots  \},
$$
and
$$
   -\alpha + \frac{1}{2} \frac{\Gamma (\alpha + 1)}{ \Gamma (\alpha -1)}
   -\frac{1}{6} \frac{\Gamma (\alpha + 1)}{\Gamma (\alpha -2)} 
   + \frac{1}{24} \frac{\Gamma (\alpha + 1)}{\Gamma (\alpha -3)}
   + \sum_{k = 3}^{\infty} w_k
   = -1,
$$
hence
\begin{displaymath}
	W_{-} = \sum_{k \in I_{-}} w_k
	= 
	   -1 - \frac{1}{2} \frac{\Gamma (\alpha + 1)}{ \Gamma (\alpha -1)}
	   - \frac{1}{24} \frac{\Gamma (\alpha + 1)}{\Gamma (\alpha -3)}
	   < 0
\end{displaymath}
and
\begin{equation} \label{eq:4.17}
   \mathbb{P} (Y^{-} = 1) =
   \frac{\alpha}
          { 1 + \frac{1}{2} \frac{\Gamma (\alpha + 1)}{ \Gamma (\alpha -1)}
	   + \frac{1}{24} \frac{\Gamma (\alpha + 1)}{\Gamma (\alpha -3)}}
\end{equation}
\begin{equation} \label{eq:4.18}
   \mathbb{P} (Y^{-} = 3) =
   \frac{  \frac{1}{6} \frac{\Gamma (\alpha + 1)}{ \Gamma (\alpha -2)} }
          { 1 + \frac{1}{2} \frac{\Gamma (\alpha + 1)}{ \Gamma (\alpha -1)}
	   + \frac{1}{24} \frac{\Gamma (\alpha + 1)}{\Gamma (\alpha -3)}}
\end{equation}
\begin{equation} \label{eq:4.19}
   \mathbb{P} (Y^{-} = k) =
   \frac{  \frac{ (-1)^{k+1} \Gamma (\alpha + 1)} {k! \, \Gamma (\alpha - k + 1)}   }
          { 1 + \frac{1}{2} \frac{\Gamma (\alpha + 1)}{ \Gamma (\alpha -1)}
	   + \frac{1}{24} \frac{\Gamma (\alpha + 1)}{\Gamma (\alpha -3)}}, 
	   \quad
	   k = 5, 6, \ldots
\end{equation}

 We arrive at the following Monte Carlo procedure.
 Let $Y_1^{+}$, $Y_2^{+}$, \ldots , $Y_n^{+}$, \ldots 
 are independent copies of random variable $Y^{+}$
 given by (\ref{eq:4.15})--(\ref{eq:4.16}),
 while $Y_1^{-}$, $Y_2^{-}$, \ldots , $Y_n^{-}$, \ldots 
  are independent copies of random variable $Y^{-}$
  given by (\ref{eq:4.17})--(\ref{eq:4.19}),
  then with probability one as $N \rightarrow \infty$
  \begin{eqnarray*}
  	& ~ &
  	\frac{1}{h^{\alpha}}
	\left[
		f(t) 
		+ \left(
			\frac{1}{2} \frac{\Gamma (\alpha + 1)}{ \Gamma (\alpha -1)}
	                 + \frac{1}{24} \frac{\Gamma (\alpha + 1)}{\Gamma (\alpha -3)}
		   \right)	
		   \frac{1}{N}
		   \sum_{n=1}^{N} f( t - Y_n^{+} h)
	\right.
		\\  
  	& ~ &
	\left.		
		- 
		\left( 
			1 +\frac{1}{2} \frac{\Gamma (\alpha + 1)}{ \Gamma (\alpha -1)}
			+ \frac{1}{24} \frac{\Gamma (\alpha + 1)}{\Gamma (\alpha -3)}
		\right)
		\frac{1}{N}
		\sum_{n=1}^{N} f( t - Y_n^{-} h)
	\right]
	\longrightarrow
	A_h^{f(t)}. 
\end{eqnarray*}

 \section{The Monte Carlo method for numerical fractional integration using continuos Sibuya-like distibution} \label{sec:6}

 Since fractional-order differentiation and integration are closely related and can be considered the unifying generalization of repeated integration and differentiation, we can consider also numerical evaluation of the Riemann-Liouville fractional integrals. 
 
 \subsection{The Riemann-Liouville fractional integration and the continuous Sibuya distribution}
 
 Let us consider the Riemann-Liouville fractional integral of order $\alpha$:

\begin{equation} \label{eq:RI-integral}
	 _{0}I_{t}^{\alpha} f(t) = \, _{0}D_{t}^{-\alpha} f(t)=
	\frac{1}{\Gamma (\alpha)} 
	\int\limits_{0}^{t} f(\tau) (t-\tau)^{\alpha-1} d\tau 
\end{equation}

Changing the variable $\tau=tu$  allows to write Riemann-Liouville fractional integral
(\ref{eq:RI-integral}) takes on the form of the Stieltjes integral: 

\begin{equation}
	_{0}I_{t}^{\alpha} f(t) = 
	\frac{t^\alpha}{\Gamma (\alpha)} 
	\int\limits_{0}^{1} f(t u) (1-u)^{\alpha-1} du =
	\frac{t^\alpha}{\alpha \Gamma (\alpha)} 
	\int\limits_{0}^{1} f(t u)  \, 
	d\left\{ 1- (1-u)^{\alpha} \right\}
\end{equation}
or
\begin{equation}
_{0}I_{t}^{\alpha} f(t) = 
\frac{t^\alpha}{\Gamma (\alpha + 1)} 
\int\limits_{0}^{1} f(t u)  
d\bigl\{ 1- (1-u)^{\alpha} \bigr\}
 = 
\frac{t^\alpha}{\Gamma (\alpha + 1)} 
\int\limits_{0}^{1} f(t  u)  \, 
dG(u),
\end{equation}
where the function
\begin{equation}
G(u) =
\left\{
\begin{array}{ll}
0, & u \in (-\infty, 0) \\
1- (1-u)^{\alpha}, & u \in [0, 1) \\
1,  & u \in [1, \infty)
\end{array}
\right.
\end{equation}
 is the cumulative function for the distribution with probability density 
 \begin{equation}
 g(u) = 
\left\{
\begin{array}{ll}
\alpha
(1-u)^{\alpha - 1}, & x \in (0, 1) \\
0, & x \not\in (0, 1)
\end{array}
\right. 
 \end{equation}

The function $G(u)$  reminds the generating function for the discrete Sibuya distribution used in the previous sections for numerical evaluation of fractional derivatives, and perhaps the distribution defined by $g(u)$  (or $G(u)$) can be called the continuous Sibyua distribution, although a close look at the above shows that this distribution is a particular case of the beta-distribution $B(\mu, \alpha)$ for $\mu =1$, 
and $G(u)$ is a particular case of the regularized incomplete beta function: $G(u) = I_u(1, \alpha)$.

 \subsection{The Monte Carlo method for Riemann-Liouville fractional integrals}
 
We have the Monte Carlo method for fractional integration similarly to the case of fractional differentiation.

\begin{itemize}

\item[] \textbf{Step 1.} Generate  $N$ random points $X_i$  using  $G(u)$.

\item[] \textbf{Step 2.} Compute 
	\begin{displaymath}
		_{0}I_{t}^{\alpha} f(t) \approx
		\frac{t^{\alpha}}
		       {\Gamma (\alpha + 1)}
			\, \,
		\frac{1}{N} 
		\sum_{i=1}^N
		f(t \, X_i)
	\end{displaymath}
	
\item[] \textbf{Step 3.} Repeat steps 1 and 2 in $K$ trials.

\item[] \textbf{Step 4.} Compute the mean of those $K$ trials, which gives the approximate value of the fractional-order integral (\ref{eq:RI-integral}). 

\end{itemize}

\section{Examples}


In the examples provided in this section 
the results of computations using the Monte Carlo method are compared 
with the exact fractional-order derivatives and integrals of the corresponding functions. 
In each figure, the solid line represents the exact fractional-order derivative or integral,
the small dots arranged vertically are the results of the repeated trials, 
the bold dot is the evaluated value (the mean of the trials), and the short horizontal lines
correspond to confidence interval for the computed value represented by the bold dot. 

In the examples provided below the Mittag-Leffler function with two parameters is used, 
$$
	E_{\delta, \beta}(z) = \sum_{k=0}^{\infty} \frac{z^k}{\Gamma (\delta k + \beta)}, 
	\quad 
	(\delta > 0, \,\,  \beta > 0),
$$
for which we have \cite[pp.~17--18]{Podlubny-FDE}:
$$
	E_{1, 1}(z) = e^z,
	\qquad
	z^2 E_{1,3}(z) = e^z - 1 - z, 
$$
$$
	z E_{2, 2} (-z^2) = \sin(z),  
	\qquad
	E_{2, 1} (-z^2) = \cos (z). 
$$
We also have a convenient rule for its fractional-order differentiation 
and integration \cite[p.~310]{Podlubny-FDE}:
%
%
$$
	_{0}D_{t}^{\alpha} 
	\left(
		t^{\beta-1} E_{\delta, \beta}(\lambda z^\delta) 
	\right)
	= t^{\beta-\alpha -1} E_{\delta, \beta -\alpha}(\lambda z^\delta). 
$$

The Mittag-Leffler function is computed using the MATLAB function~\cite{Podlubny-MLF}, 
and the examples are included in~\cite{Podlubny-MCFD-toolbox}.

\subsection{Example 1: fractional differentiation of order  $0 < \alpha < 1$. }

\begin{equation}
	y(t) = E_{\alpha, 1} (-\lambda t^\alpha) - 1
\end{equation}

\begin{equation}
	_0D_t^{\alpha} y(t) = t^{-\alpha} E_{\alpha, 1-\alpha} (-\lambda t^\alpha) 
	                            - \frac{t^{-\alpha}}{\Gamma (1+\alpha)}
\end{equation}

\begin{figure}
\begin{center}
	\includegraphics[width=0.6\textwidth]{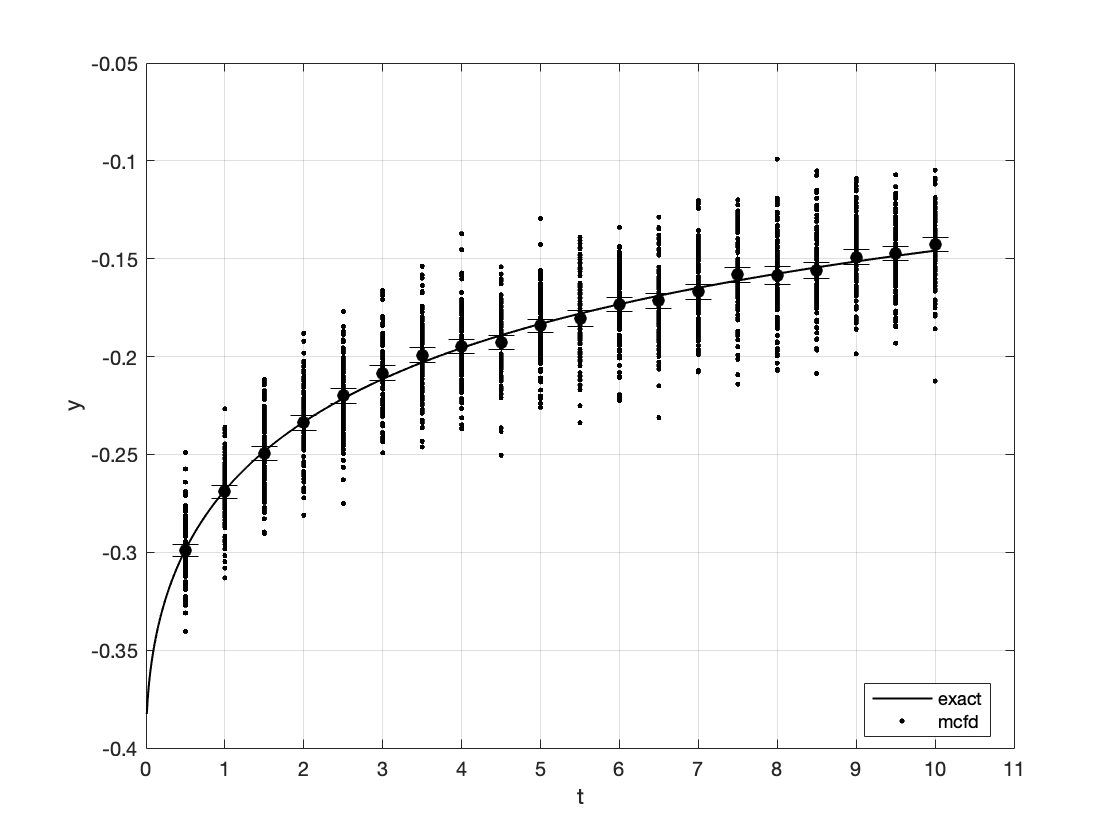}\\	
\caption{Evaluation of $_0D_t^{\alpha} y(t)$ for $y(t) = E_{\alpha, 1} (-\lambda t^\alpha) - 1$; $\alpha = 0.5$, $\lambda = 0.4$ }
\label{fig:example1}
\end{center}
\end{figure}

The result of evaluation for 
$\alpha = 0.5$, $\lambda = 0.4$
is shown in Fig.~\ref{fig:example1}.
	
\subsection{Example 2: fractional differentiation of order  $1 < \alpha < 2$. }

\begin{equation}
	y (t) = \sin (t)
\end{equation}

\begin{equation}
	_0D_t^{\alpha} y(t) = t^{1-\alpha} E_{2, 2-\alpha} (-t^2)
\end{equation}

\begin{figure}
\begin{center}
	\includegraphics[width=0.6\textwidth]{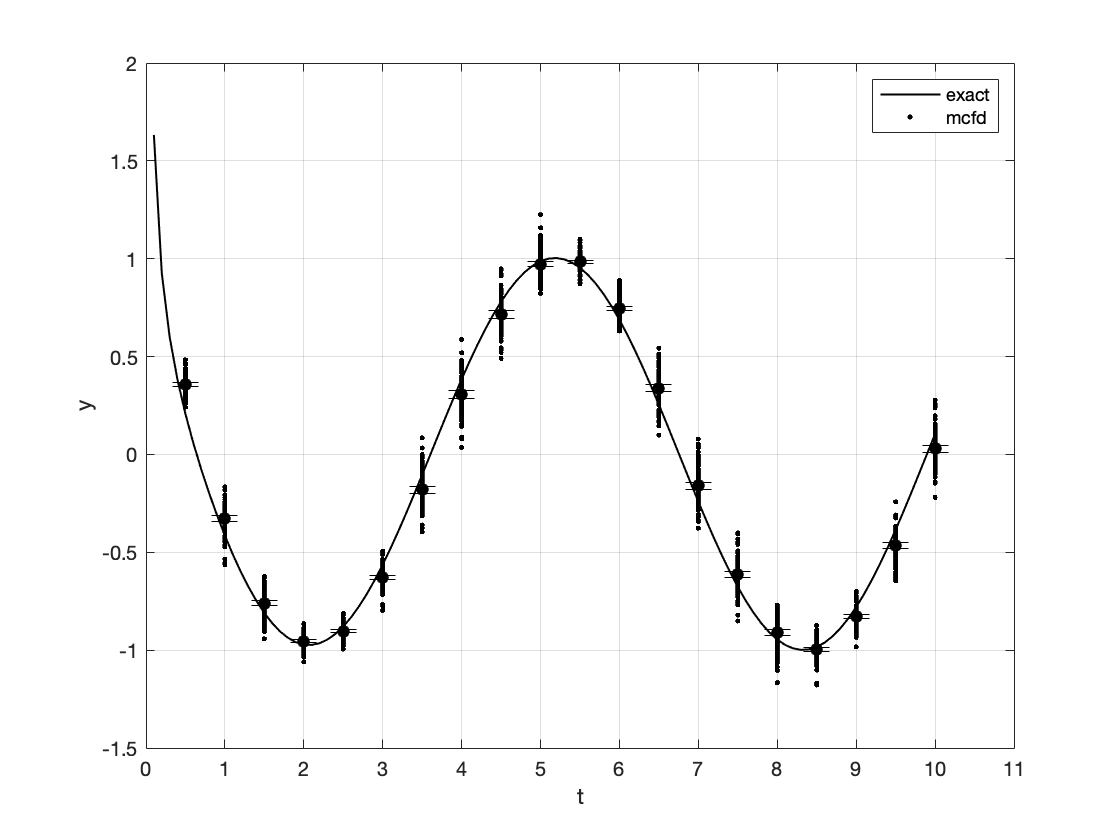}\\	
\caption{Evaluation of $_0D_t^{\alpha} y(t)$ for $y(t) = \sin (t)$; $\alpha = 1.7$ }
\label{fig:example2}
\end{center}
\end{figure}

The result of evaluation for 
$\alpha = 1.7$
is shown in Fig.~\ref{fig:example2}.

\subsection{Example 3: fractional differentiation of order  $2 < \alpha < 3$.}

\begin{equation}
	y(t) = e^{-\lambda t} - 1 + \lambda t; 
\end{equation}

\begin{equation}
	_0D_t^{\alpha} y(t)  = \lambda^2 t^{2-\alpha} E_{1, 3-\alpha} (-\lambda t)
\end{equation}

\begin{figure}
\begin{center}
	\includegraphics[width=0.6\textwidth]{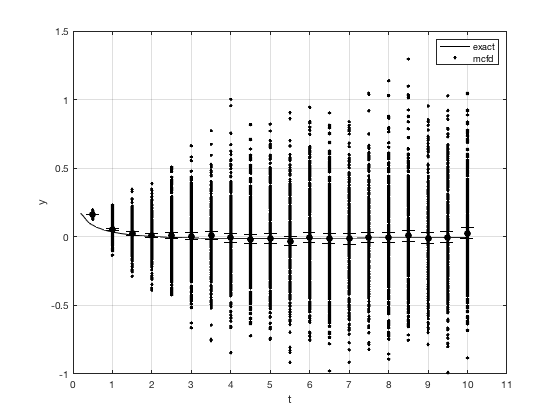}\\	
\caption{Evaluation of $_0D_t^{\alpha} y(t)$ for $y(t) = e^{-\lambda t} - 1 + \lambda t$; $\alpha = 2.5$ }
\label{fig:example3}
\end{center}
\end{figure}

The result of evaluation for 
$\alpha = 2.5$
is shown in Fig.~\ref{fig:example3}.

\subsection{Example 4: fractional integration of order  $0 < \alpha < 1$.}

\begin{equation}
	y(t) = t^{-\alpha} E_{2, 1-\alpha} (-t^2)
\end{equation}

\begin{equation}
	_{0}I_{t}^{\alpha} y(t) = \cos(t)
\end{equation}

\begin{figure}
\begin{center}
	\includegraphics[width=0.6\textwidth]{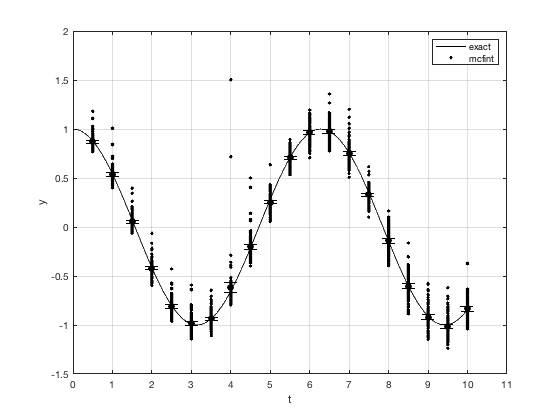}\\	
\caption{Evaluation of $_0I_t^{\alpha} y(t)$ for $y(t) = t^{-\alpha} E_{2, 1-\alpha} (-t^2)$; $\alpha = 0.6$ }
\label{fig:example4}
\end{center}
\end{figure}

The result of evaluation for 
$\alpha = 0.6$
is shown in Fig.~\ref{fig:example4}.

\subsection{Example 5: fractional integration of order  $1 < \alpha < 2$.}

\begin{equation}
	y(t) = t^{1-\alpha} E_{2, 2-\alpha} (-t^2)
\end{equation}

\begin{equation}
	_{0}I_{t}^{\alpha} y(t) = \sin(t)
\end{equation}

\begin{figure}
\begin{center}
	\includegraphics[width=0.6\textwidth]{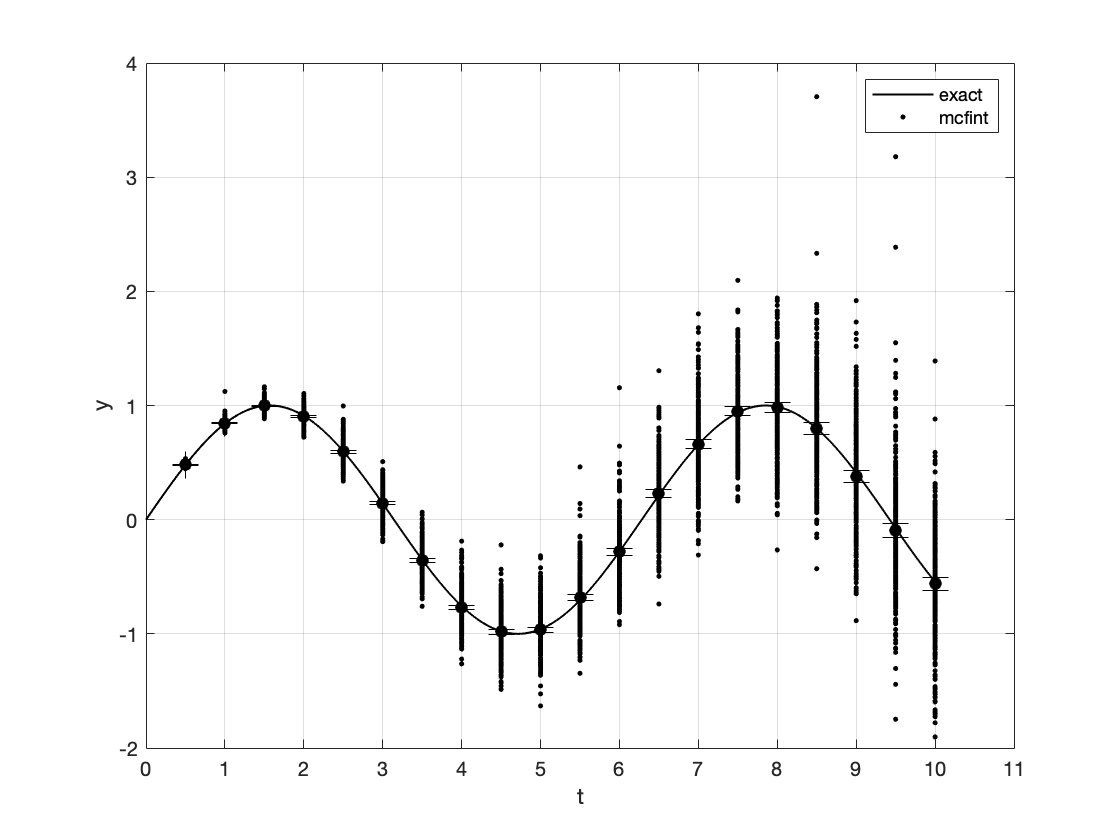}\\	
\caption{Evaluation of $_0I_t^{\alpha} y(t)$ for $y(t) = t^{1-\alpha} E_{2, 2-\alpha} (-t^2)$; $\alpha = 1.4$ }
\label{fig:example5}
\end{center}
\end{figure}

The result of evaluation for 
$\alpha = 1.4$
is shown in Fig.~\ref{fig:example5}.

\subsection{Example 6: fractional integration of order  $2 < \alpha < 3$.}

\begin{equation}
	y(t) = \frac{t^{\nu}}{\Gamma (\nu + 1)} 
\end{equation}

\begin{equation}
	y(t) = \frac{t^{\nu + \alpha}}{\Gamma (\nu +  \alpha + 1)} 
\end{equation}

\begin{figure}
\begin{center}
	\includegraphics[width=0.6\textwidth]{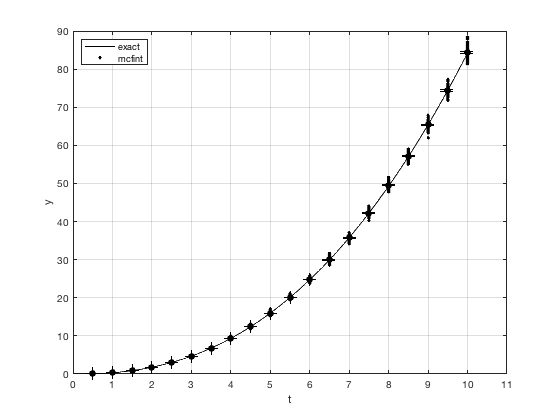}\\	
\caption{Evaluation of $_0I_t^{\alpha} y(t)$ for $y(t) = \frac{t^{\nu}}{\Gamma (\nu + 1)}$; $\nu = -0.3$, $\alpha = 2.7$ }
\label{fig:example6}
\end{center}
\end{figure}

The result of evaluation for 
$\nu = -0.3$ and $\alpha = 2.7$
is shown in Fig.~\ref{fig:example6}.

\section{Conclusion} 

We can conclude that, due to its link to fractional-order differentiation and integration, 
the Sibuya distribution deserves close attention and investigation. 
Consideration of the case when its parameter is greater than one 
takes us to the field of signed probability distribution. 

Because of the link between fractional-order differentiation and fractional-order integration,
we can think of the convenient unification of tools and terminology, 
and attribute the name of the continuous Sibuya distribution
to a particular case of the beta distribution. 

For higher orders of fractional differentiation, the simulation of the sieved Sibuya distributions brings computational challenges. Also, the presented Monte Carlo method using signed probabilities can be further enhanced with the help of standard approaches 
for improving the classical Monte Carlo method. Both aspects require further exploration.

\end{document}